\documentclass[authoryear]{elsarticle}
\RequirePackage[english]{babel}
\RequirePackage{lmodern}
\RequirePackage{amsmath}
\RequirePackage{amssymb}
\RequirePackage{amsthm}
\RequirePackage{amsfonts}
\RequirePackage{algorithm}
\RequirePackage[noend]{algorithmic}
\RequirePackage{enumerate}
\RequirePackage{enumitem}
\RequirePackage{mathrsfs}
\RequirePackage{multirow}
\RequirePackage{multicol}
\RequirePackage{nicematrix}
\RequirePackage{pifont} 
\RequirePackage{rotating}
\RequirePackage{tikz}
\usetikzlibrary{arrows.meta}
\usepackage[normalem]{ulem} 

\RequirePackage[pdftex,pdfpagelabels,bookmarks,hyperindex,hyperfigures]{hyperref}
\hypersetup{
	colorlinks=true,
	linkcolor=blue,
	citecolor=blue,
}

\newtheorem{example}{Example}
\newtheorem{property}{Property}
\newtheorem{theorem}{Theorem}
\newtheorem{definition}{Definition}
\newtheorem{assumption}{Assumption}
\newtheorem{lemma}{Lemma}
\newtheorem{corollary}{Corollary}
\newtheorem{remark}{Observation}

\def\squareforqed{\hbox{\rlap{$\sqcap$}$\sqcup$}}
\def\qed{\ifmmode\squareforqed\else{\unskip\nobreak\hfil
		\penalty50\hskip1em\null\nobreak\hfil\squareforqed
		\parfillskip=0pt\finalhyphendemerits=0\endgraf}\fi}

\newcommand{\constraint}[1]{\textsc{#1}\,}

\newcommand{\ie}{\textit{i.e. }}

\newcommand{\st}{\textrm{s.t. }}

\newcommand{\defined}[1]{\textbf{\boldmath#1}}

\newcommand{\set}[1]{\ensuremath{\mathcal{#1}}}

\newcommand{\Z}{\ensuremath{\mathbb{Z}}}

\newcommand{\dual}{\ensuremath{\mathcal{D}}}
\newcommand{\primal}{\ensuremath{\mathcal{P}}}

\newcounter{constraint}
\newcommand{\prooflabel}[1]{\hypertarget{#1}{}}

\tikzset{
	node/.style={circle, draw, fill=lightgray},
	arc/.style={-{Latex[length=2.5mm]}}
		endnode/.style={line width=1pt},
		matching/.style={red},
		nrctree/.style={line width=2pt},
		filtered/.style={line width=1.5pt, dashed}
}

\title{Arc-consistency with linear programming reduced costs\\ (applied to stable set in chordal graphs)}

\author[1]{Guillaume Claus}\ead{}
\author[1]{Hadrien Cambazard}
\ead{}
\author[2]{Hugo Apeloig} \ead{}
\author[1]{Pierre Hoppenot}
 \ead{hadrien.cambazard@grenoble-inp.fr}
\ead{}

\affiliation[1]{organization={Univ.Grenoble Alpes, CNRS, Grenoble INP, G-SCOP}, city={Grenoble}, country={France}}
\affiliation[2]{organization={Université de Nantes, IMT Atlantique}, city={Nantes}, country={France}}

\begin{document}
\begin{abstract}

	A well known technique to reduce the search space in integer programming is known as \emph{variable fixing} or \emph{reduced cost strengthening}. The reduced costs given by an optimal dual solution of the linear relaxation can be used to strengthen the bounds of the variables but this filtering is incomplete. We show how reduced costs can be used to achieve Arc-Consistency (AC), \emph{i.e} a complete filtering, of a global constraint with a cost variable and an assignment cost for each value. We assume that an ideal Integer Linear Programming (ILP) formulation is available i.e the convex hull of the characteristic vectors of the supports is known. A detailed analysis of reduced cost based filtering is proposed. We characterize arc-consistency based on complementary slackness \emph{i.e} completeness of reasoning as opposed to only optimality. We also give a simple sufficient condition allowing a set of dual solutions to ensure arc-consistency through reduced costs.
	In practice, when the constraint has a such an ideal ILP, $n$ dual solutions are always enough to achieve AC (where $n$ is the number of variables of the global constraint). It extends the work presented in \citep{german2017arc} for satisfaction problems and in \citep{claus2020wad} for the specific case of the minimum weighted alldifferent constraint.
	Our analysis is illustrated on constraints related to the assignment and shortest path problem and also demonstrated on the weighted stable set problem in chordal graphs. A novel AC algorithm is proposed in this latter case based on reduced costs. 
\end{abstract}
\maketitle

\section{Introduction}	
	Mixed Integer Programming (MIP) and Constraint Programming (CP) have often been combined in the past to take advantage of the complementary strengths of the two frameworks.  Many approaches have been proposed to benefit from their modelling and solving capabilities~\citep{Bockmayr1998,Rodosek1999,refalo2000linear,aron2004simpl,achterberg2008constraint}. A typical integration of the two  approaches is to use the linear relaxation of the entire problem in addition to the local consistencies enforced by the CP solver.  The relaxation can detect infeasibility and is often added to provide a bound on the objective.
	
	A number of previous works have also proposed to use the linear relaxation for filtering the domains in a constraint programming framework \citep{refalo1999tight,refalo2000linear,aron2004simpl,achterberg2008constraint,FocacciLM02}. Based on the relaxation, filtering can be performed using a technique referred to as reduced cost based filtering \citep{FocacciLM02,HookerMS06}. It is a specific case of cost based filtering \citep{focacci1999cost} that aims at filtering out values leading to non-improving solutions. It originates from \emph{variable fixing} \citep{Nemhauser:1988:ICO:42805} which is performed in MIP to detect some 0/1 variables that must be fixed to either 0 or 1 in any solution improving the best known. Variable fixing usually relies on the reduced costs of the variables given by an optimal dual solution of the linear relaxation. It was already used in 1954 (Section \emph{An estimation procedure} of \citep{dantzig1954solution}) for solving the Traveling Salesman Problem (TSP) and early in integer programming \citep{balas1980pivot,crowder1983solving} as well as Lagrangian relaxation \citep{beasley1990lagrangian}. Reduced cost based filtering is a typical component of modern MIP solvers as one of the bound tightening techniques \citep{gleixner2017three}. It is known to be incomplete because it strongly depends on the specific dual solution used and many authors are investigating how to search the dual space to gather as much filtering as possible from the reduced costs \citep{sellmann2004theoretical}. To do so, a dual picking program is proposed in \citep{bajgiran2017first}. In the context of Lagrangian relaxation, \citep{boudreaultimproved, DBLP:conf/cpaior/BerthiaumeQ24} design an algorithm that alters locally the Lagrange multipliers to enhance reduced cost filtering by focusing on promising variables. A recent application to network flows \citep{de2023exact} is using three heuristics to produce dual values increasing the filtering. It is now-days widely used in practice in many areas from vehicle routing \citep{schurmann2023reduced} to maxSAT \citep{bacchus2017reduced}. Alternatively, it was shown in \citep{german2017arc} that a complete filtering, namely arc-consistency, can be achieved by solving a single linear relaxation when the problem considered is a satisfaction problem with a known ideal Integer Linear Programming (ILP) formulation. More precisely, this ILP is based on 0/1 variables that typically encode whether an integer variable of the constraint’s scope is assigned to a value of its domain. It is required to be \emph{ideal} \emph{i.e} that its linear relaxation defines a polytope representing the convex hull of the characteristic vectors of the supports of the constraint. Such formulations can be found for a number of common global constraints such as \constraint{Element}, \constraint{AllDifferent}, \constraint{GlobalCardinality} or \constraint{Gen-Sequence} \citep{refalo2000linear,german2017arc}.
	The approach does not apply to global constraints involving a cost variable such as \constraint{MinimumWeightAlldifferent} \citep{DBLP:journals/constraints/CaseauL00,focacci1999cost} even though a compact ideal formulation of the matchings polytope is available. A natural extension to the work \citep{german2017arc} is to handle an objective function \ie a cost variable from the constraint point of view. We are therefore interested in the design of filtering algorithms based on linear programming for global constraints with a cost variable. Note that when an ideal ILP is available for the constraint, a naive approach, typically used in practice when debugging and prototyping propagators is to solve one LP for each variable-value pair. Such a technique also relates to shaving initially introduced in a scheduling context \citep{CARLIER1994146}, or probing \citep{savelsbergh1994preprocessing}, in particular when it is done for the entire problem in CP and not on a single constraint.

	We consider global constraints with assignment costs. More precisely, assigning a value to a variable incurs a cost and the overall cost is the sum of all individual costs. Soft global constraints might have alternative costs definition but assignment costs are very common.
    The consistency of a given value of a variable can be established by computing the minimum objective value of the problem restricted to the corresponding assignment. When this increase of the objective is inconsistent with the upper bound of the cost variable, the value is inconsistent. A reduced cost is a typical lower bound of this increase and is obtained for free from a single optimal dual solution as a by-product. This filtering is computationally cheap and was used in CP to perform filtering in \citep{focacci1999cost} for the assignment problem, but it remains incomplete. The reduced cost depends on the optimal dual solution found so when it is not unique, reduced costs tend to greatly vary in practice from one solution to another. When a reduced cost gives the exact increase of the objective, we refer to it as exact.\\

	The contribution of this paper is twofold. Firstly, we provide a novel analysis of reduced cost based filtering by relating arc-consistency to complementary slackness and linear programming duality theory. This is illustrated on the assignment and shortest path problem in an acyclic graph. We consider constraints for which an ideal ILP formulation is known. 
	We prove (Property \ref{prop: exists u st r_kl=R_kl}) that for a given variable/value pair, there always exists an optimal dual solution such that its reduced cost is exact. We give necessary and sufficient conditions (Theorem \ref{ac slackness} and Corollary \ref{prop: characterization r_kl=R_kl}) based on complementary slackness (CS) to identify which reduced costs are exact in a given optimal dual solution. The relationship with CS eventually leads to a sharp upper bound on the number of dual solutions needed to ensure AC (Property \ref{prop:upperBoundNbDSol}). 
	Given a set of variable-value pairs, we give a sufficient condition (Property \ref{optSolDtilde} and Corollary \ref{prop: incompatible edges}) for the existence of a single optimal dual solution providing all exact reduced costs for this set. More precisely, exact reduced costs of values that are pairwise inconsistent \emph{i.e} that do not belong simultaneously to a feasible solution, can be found in the same dual solution. This condition leads to an algorithm based on LP to enforce AC (Algorithm \ref{algoACbyLP}). This generic approach applies to all global constraints with assignment costs and such an ideal ILP. 
	
	The previous work of \citep{german2017arc} was based on an interior point of the primal formulation. In the present work, this result is presented as an interior point of the dual (Property \ref{prop: r_ij>0 <=> ij filterable}) relating the analysis to reduced costs which are the slacks of the dual constraints. 
	Moreover, we now consider weighted constraints with a cost variable adding an objective function to the LP for identifying an optimal support. This extension requires several dual solutions that must be carefully chosen to ensure AC. Compared to \citep{claus2020wad}, the presented results are more general than the case of \constraint{MinimumWeigthAllDifferent} and, more importantly, related to complementary slackness conditions leading to a new characterization of an exact reduced cost.  To be precise, Properties \ref{optSolDtilde} to \ref{prop:upperBoundNbDSol} were presented in \citep{claus2020wad} for a specific case and now generalized. Theorems \ref{ac slackness}, \ref{simultaneous exact rc} and Corollary \ref{prop: characterization r_kl=R_kl} are novel.
	
	Secondly, we apply these results to design a complete filtering algorithm for the weighted maximum stable set problem in chordal graphs. This algorithm is novel and demonstrates the interest of the previous statements for a fundamental problem. It shows how a combinatorial algorithm relying on duality can be adapted to enforce arc-consistency.\\


	In Section \ref{section: AC for global constraint with costs} we set the framework for our work and present the two constraints used as examples:  \constraint{MinimumWeight\-Alldifferent} (referred to as \constraint{MinWAllDiff} for short in the rest of the paper), and \constraint{ShortestPath}.	
	The LPs  used in this document and the relationship between filtering and reduced costs are explained in Section \ref{section: LP formulations}. 
	The main results are detailed at Section \ref{section: analysis} and the resulting filtering algorithm is stated in Section \ref{section: LP base algo}. Finally, Section \ref{section: MWIS} demonstrates an application of these results on a global constraint \constraint{MaxWIS} dealing with Maximum Weighted Independent Sets in chordal graphs.
	
\section{Preliminaries on Arc-Consistency}
\label{section: AC for global constraint with costs}

	A \defined{constraint satisfaction problem} (CSP) is made of a set of variables, each with a given \defined{domain} \ie a finite set of possible values, and a set of constraints specifying the allowed combinations of values for subsets of variables. 
	In the following, the variables, \emph{e.g.} $\{X_i\}_{i=1}^n$, are written with upper case letters for the constraint programming models as opposed to the variables of linear programming models that are in lower case.  
	$D(X_i) \subset \Z$ denotes the domain of $X_i$.
	The minimum and maximum values in $D(X_i)$ are respectively denoted $\underline{X_i}$ and $\overline{X_i}$. 
	A \defined{constraint} $C$ over a set of variables $\langle X_1, \ldots, X_n \rangle$ is defined by the allowed combinations of values (tuples) of its variables. Such tuples of values are also referred to as solutions of the constraint $C$.
	Given a constraint $C$ with a scope $\langle X_1, \ldots, X_n \rangle$, a \defined{support} for $C$ is a tuple of values $\langle a_1,\ldots, a_n \rangle$ that is a feasible solution of $C$ and such that $a_i \in D(X_i)$ for all variables $X_i$ in the scope of $C$. 
	Consider a variable $X_i$ in the scope of $C$, the domain $D(X_i)$ is \defined{arc-consistent} for $C$ if and only if all the values of $D(X_i)$ belong to a support for $C$. A constraint $C$ is  arc-consistent if and only if all its variable's domains are arc-consistent. 
	
	Let $C\left(X_1, \dots ,X_n\right)$ be a constraint. Given a cost $c_{ij} \in \mathbb{Z}$ for assigning variable $X_i$ to value $j \in D(X_i)$, we define $W\!C\left(X_1, \dots ,X_n,Z,c\right)$ as $C(X_1,\ldots,X_n)\:\wedge\:\sum^{n}_{i=1} c_{i,X_i} \leq Z$ where $Z$ is a cost variable ($D(Z) \subset \mathbb{Z}$). We refer to $W\!C$ as a \textbf{weighted constraint} and its associated constraint $C$ as the constraint without cost.
	The cost of a tuple $\langle a_1,\ldots, a_n \rangle$ is defined as $\sum^{n}_{i=1} c_{i,a_i}$. 
	The constraint holds if $\langle a_1,\ldots, a_n \rangle$ is a support for $C$, and its cost is below $\overline{Z}$. For sake of simplicity, to state the general results, we consider weighted constraints enforcing an upper bound of the cost. In other words, $\langle a_1,\ldots, a_n, z \rangle$ is a \defined{support for $W\!C$} if $\langle a_1,\ldots, a_n \rangle$ is a support for $C$ and $z = \sum_{i=1}^n c_{i,a_i} \leq \overline{Z}$. Typically, the filtering algorithm of $W\!C$ updates $\underline{Z}$ by solving a minimization problem and filters values that are inconsistent for $C$ or too costly regarding $\overline{Z}$.
	
	A \defined{support of minimal cost} is a support for $W\!C$ such that its cost is no more than the cost of any other support for $W\!C$. It is a support for $\underline{Z}$ and its cost is denoted $z^*$.
	\begin{example}
		\slshape\sffamily
		~
		 \noindent\constraint{MinWAllDiff}$(X_1,\ldots, X_n, Z, c)$ is equivalent (has the same set of solutions) to the constraint network \citep{sellmann2002arc}:
			$$\begin{array}{ll}
				\constraint{AllDiff}(X_1,\ldots, X_n) \\
				\sum^{n}\limits_{i=1} c_{i,X_i} \leq Z \\
				X_i \in D(X_i) \qquad \forall i \in \{1,\ldots,n\}
			\end{array}$$
		
			\noindent Where \constraint{AllDiff} ensures the $X$ variables take distinct values.

			\indent\constraint{ShortestPath}$(X_1,\dots,X_{n},Z,c)$ is equivalent to the constraint network:
			$$\begin{array}{l}
				\constraint{Path}(X_1,\dots,X_{n})\\
				\sum\limits_{i=1}^n c_{i,X_i}  \leq Z\\
				X_i \in \{i\} \cup \{j | ij \in \mathcal{A}\} \:\:\: \forall i \in \mathcal{N}\setminus\{t\}
			\end{array}
			\textrm{with } X_{i} =\begin{cases}
				j &  \text{if $j$ is the successor of $i$}\\
				i & \text{if $i$ is not in the path}
			\end{cases}$$
			\noindent Note that $c_{ii}=0$ for all $1\le i \le n$. \constraint{Path} ensures the $X$ variables describe a path from vertex $1$ to vertex $n+1$ in an acyclic graph denoted $(\mathcal{N},\mathcal{A})$. Vertices $1$ and  $n+1$ act as the source $s$ and sink $t$ ($s=1$ and $t = n+1$). So $\mathcal{N}$ is made of $n+1$ vertices, one for each variable and the sink. The set of arcs correspond to the set of values of the domains ignoring the loops $ii$ \emph{i.e} $ij \in \mathcal{A} \Leftrightarrow j \in D(X_i), i \neq j$
	
		Another possible formulation for \constraint{ShortestPath} constraint is to use 0/1 variables, one for each possible arc. Each of these variables equals to 1 if and only if the path goes through the corresponding arc. In this case, dummy value $i$ for $X_i$ is not encoded.
	\end{example}
We denote by $ij$ (short for $(i,j)$ to simplify notations) a \emph{variable/value pair} $(X_i, j)$ such that value $j \in D(X_i)$. The set $B$ is the set of all possible pairs given by the variables $U = \{X_1, \ldots, X_n\}$ in the scope of $C$. So $B = \{ij \: | \: \forall i \in U, \forall j \in D(X_i)\}$. We refer to a \textbf{pair} $ij$ but also, by abuse of language, to a \textbf{value} $ij$. The cost of a minimal support that includes the pair $ij$ (such that $X_i = j$) is denoted $z^*_{|ij}$. 


The problem of identifying a minimum support of $W\!C$ can be stated in integer linear programming and such a formulation is referred to as a \textbf{formulation of $W\!C$}.

The \textbf{characteristic vector} of a support is a 0/1 vector of $\{0,1\}^{|B|}$ encoding whether a pair $ij$ belongs or not to the support.  Considering a weighted constraint $W\!C$, we denote by $\mathcal{S} \subseteq \{0,1\}^{|B|}$ the set of the characteristics vectors of the supports of $C$. A formulation of $W\!C$ can be stated as $\min\{cx : x\in \mathcal{S}\}$. In the following, we assume that a linear description of the convex hull, $conv(\mathcal{S})$, of $\mathcal{S}$ is available. The formulation $\min \{cx: x \in conv(\mathcal{S}), x \in \{0,1\}^{|B|}\}$ is said ideal (\cite{Wolsey:1998} p. 15) because the integer problem can be solved by solving its linear relaxation $\min \{cx: x \in conv(\mathcal{S})\}$ since each extreme point is integer. In the following and by abuse of language, $\min \{cx: x \in conv(\mathcal{S})\}$ as well as any formulation whose feasible region is $conv(\mathcal{S})$ are also referred to as \textbf{ideal formulations}. 

 
\section{Linear programs of global constraints}
\label{section: LP formulations}
	  
Let $A$ be an $m \times |B|$-matrix and $b$, a row vector of dimension $m$. We assume an ideal formulation of the constraint $W\!C$ is available and feasible region is defined by a polytope $$L = \{x \in \mathbb{R}^{|B|}_+ \:\: : \:\: Ax \ge b\} = conv(\mathcal{S})$$
Integer points of $L$ are exactly the characteristic vectors of the supports of $C$. Note that the variables of $L$ provide an encoding of the integer domains of the constraint  where a variable $x_{ij}$ is used for each pair $ij \in B$ so that: $x_{ij} = 1 \iff X_i=j$.
Note also that the formulation ensures that the constraints of the domains enforce each variable to take a single value \emph{i.e} $\sum_{j|\,ij \in B} x_{ij} = 1\quad \forall i \in U$. 
	The problem of identifying a support of minimal cost $z^*$ for $W\!C$ is stated as :
	$$(\primal) \qquad z^* = \min \{cx : x \in L\}$$
	Since all extreme points are integral, there is an integer solution $x^*$ of value $z^*$. More generally, $z$ denotes the cost of a solution of $(\primal)$ and this solution is a support for $W\!C$ if it is integer and $z \le \overline{Z}.$ 
	In practice, when some values are functionally dependent from others, formulation $(\primal)$ can be stated with less variables than $|B|$ (see example \ref{exampleLP} with $(\primal_{\constraint{SP}})$ below).
	
		A dual of $(\primal)$ is given by:
	$$(\dual) \qquad w^* = \max \{ub : u A \leq c, u \in \mathbb{R}^{m}_+\}$$

	\begin{example}[continued: Linear Programs]
	\label{exampleLP}
		\slshape\sffamily
		
		Formulation $\left(\primal\right)$ is given below for \constraint{MinWAllDiff} and denoted $(\primal_{\constraint{WAD}})$. It is a typical formulation of the linear assignment problem on the weighted bipartite graph $(U \cup V, E, c)$ where the set $U$ refers to the set of the variables whereas $V$ is the set of values, common for all the variables \emph{i.e} $V=\cup_{i \in U} D(X_i)$. This formulation was initially used in CP by \citep{focacci1999cost} and is known to be ideal (Section 4.3 in \cite{Wolsey:1998}). A support, in this case, is a matching of cardinality $n$.
	
		\begin{center}
		\begin{tabular}{ll}	
			$(\primal_{\constraint{WAD}})$ &
			$\left\{	
			\begin{array}{rrll@{\:\:}l@{\:\:}l@{\:\:}l}
				\multicolumn{5}{l}{\min z = \sum\limits_{ij\in E} c_{ij} \, x_{ij}}\\
				\text{s.t. }
				& \sum\limits_{j|\,ij\in E} x_{ij} & = & 1 & \forall i \in U \\
				& \sum\limits_{i|\,ij \in E} x_{ij} & = & 1\quad & \forall j \in V\\
				& x_{ij} & \ge &  0 &\forall ij \in E 
			\end{array}\right.$
		\end{tabular} 
	    \end{center}
	    
		Formulation $\left(\primal\right)$ is given below for \constraint{ShortestPath} and denoted $(\primal_{\constraint{SP}})$. A support of \constraint{ShortestPath} is an $s$-$t$ path of $(\mathcal{N},\mathcal{A})$ extended with the values $i$ of $X_i$ for all vertices $i$ that do not belong to the path. The pairs $ii$ in a characteristic vector of a support are known once the path is known and we focus on the characteristic vectors of the paths alone. The consistency of a value $i$ for $X_i$ can then be derived once the consistency of all values $j \in D(X_i)$ with $j \neq i$ is established. We do not discuss this in details since it is not the focus of the example. It is possible under Assumption \ref{feasAssump} (presented below) stating that any value in the domain belongs to a feasible support (independently of its cost) i.e any arc belong to at least one s-t path visiting the mandatory vertices ($i$ s.t $i \not\in D(X_i)$).  Under this assumption, we can focus on the characteristic vectors of the paths alone (ignoring the pairs $ii$).

		The formulation $(\primal_{\constraint{SP}})$ is a typical flow formulation of the shortest path problem stated here in the acyclic graph $(\mathcal{N},\mathcal{A})$. Variables $x_{ii}$ encoding the dummy value of each $X_i$ are not explicitly included in  $(\primal_{\constraint{SP}})$ and not required. $(\primal_{\constraint{SP}})$ is known to be ideal (Section 3.4.1 in \cite{Wolsey:1998}).
		
		\begin{center}
		\begin{tabular}{ll}
			$(\primal_{\constraint{SP}})$  & $\left\{	
			\begin{array}{rrll@{\:\:}l@{\:\:}l@{\:\:}l}
				\multicolumn{5}{l}{\min z = \sum\limits_{ij\in A} c_{ij} \, x_{ij}}\\
				\text{s.t. } 	
				& \sum\limits_{j|\,kj\in A} x_{kj} - \sum\limits_{i|\,ik \in A} x_{ik} & = & 0 & \forall k \in \mathcal{N}\setminus\{s,t\}\\
				& \sum\limits_{j|\,sj \in A} x_{sj} & = & 1\\
				& x_{ij} & \ge &  0 &\forall ij \in \mathcal{A} 
			\end{array}\right.$ %
		\end{tabular}
		\end{center}
		A support of minimal cost is obtained from an optimal solution to $(\primal_{\constraint{SP}})$ as follows. We set $X_i = j$ whenever $x_{ij} = 1$ and $X_i = i$ when $x_{ij} = 0 \:\: \forall j \textrm{ s.t } ij \in \mathcal{A}$.
		
		Both these formulations are known to be ideal since their constraint matrix ($A$) is known to be totally unimodular, the right hand side of the constraints are integer values and the lower bounds of the variables are also integers (equal to 0).
	\end{example}
	
	\medskip

	\begin{example}[continued: dual linear programs]
		\slshape\sffamily~
		
		Dual formulations are given below for  \constraint{MinWAllDiff} and \constraint{ShortestPath}.\\
		
		\noindent$\begin{array}{l|l}	 
		(\dual_{\constraint{WAD}}) \qquad \qquad \qquad \qquad \qquad \qquad \qquad & (\dual_{\constraint{SP}}) \\
				\left\{				
			\begin{array}{rrll@{\quad}l@{\qquad}l@{\qquad}l}
				\multicolumn{6}{l}{\max w = \sum\limits_{i\in U}  u^1_i + \sum\limits_{j\in V} u^2_j}\\
				\text{s.t. } 	
				& u^1_i + u^2_j & \le & c_{ij} & \forall ij\in E\\
				& u^1_{i} & \in &  \mathbb{R} & \forall i \in U\\
				& u^2_{j} & \in &  \mathbb{R} & \forall j \in V
			\end{array}\right.
			&
			\left\{				
			\begin{array}{rrll@{\quad}l@{\qquad}l@{\qquad}l}
				\multicolumn{6}{l}{\max w = u_s}\\
				\text{s.t. } 	
				& u_i - u_j & \le & c_{ij} & \forall ij\in \mathcal{A}, j \neq t\\
				& u_i & \le &  c_{it} & \forall it \in \mathcal{A}\\
				& u_i & \in &  \mathbb{R} & \forall i \in \mathcal{N}\setminus\{t\}
			\end{array}\right.
		\end{array}$
	\end{example}
	Additionally, let $\left(\primal_{|kl}\right)$ be the LP defined by $\left(\primal\right)$ in which $x_{kl}$ is forced to 1 (\ie $\left(\primal_{|kl}\right)$ is the restricted problem with the additional constraint $x_{kl}=1$), and $z_{|kl}^*$ its optimal value. Let us assume $\mathcal{S}_{|kl}$, the set of supports where the pair $kl$ is set to 1, is not empty (stated more precisely with Assumption \ref{feasAssump} below). Formulation $(\primal_{|kl})$ is an ideal formulation since its polytope $L_{|kl}$  defined as $L_{|kl} = \{x \in \mathbb{R}^{|B|}_+ \:\: : \:\: Ax \ge b, x_{kl} = 1\}$ is the convex hull of the supports $\mathcal{S}_{|kl}$ so that $L_{|kl} = conv(\mathcal{S}_{|kl})$. 
	
Using the concepts of \cite{schrijver1998theory} p.99-111, $L_{|kl} = conv(\mathcal{S}_{|kl})$ can be justified as follows. Consider $H_{|kl} = \{x\in \mathbb{R}^{|B|} : x_{kl}=1\}$. We have  $\mathcal{S}_{|kl} = \mathcal{S}\cap H_{|kl}$ not empty and $L_{|kl} = L \cap H_{|kl}$. Since $H_{|kl}$ is a supporting hyperplane of $L$, $L_{|kl}$ is a face of $L$ and $L_{|kl} = conv(\mathcal{S}_{|kl})$. 
	
	Finally, $\left(\dual_{|kl}\right)$ refers to the dual of $\left(\primal_{|kl}\right)$. 
	For the rest of this document, primal and dual solutions refer to solutions of $\left(\primal\right)$ and $\left(\dual\right)$ when it is not specified otherwise. 
	
	\begin{assumption}[Feasible and bounded formulations] \label{feasAssump}  Consider a weighted constraint $W\!C$ and assume that AC has been achieved on the constraint without costs $C\left(X_1,\dots, X_n\right)$ related to $W\!C$. In other words, any pair belongs to at least one support of $C\left(X_1,\dots, X_n\right)$, ignoring the costs. So $(\primal)$ is assumed feasible as well as all $(\primal_{|kl})$ for all pairs $kl \in B$. Note that $(\primal)$ and all $(\primal_{|kl})$ are also bounded LPs because the initial domains are finite and initial costs $c_{ij}$ are finite. So $z^*$ as well as $z^*_{|kl}$ are finite by assumption. As a result, $(\dual)$ and all $(\dual_{|kl})$ are also feasible and bounded LPs.
	\end{assumption}


	The \defined{reduced cost} of a pair $kl \in B$ is the slack of the corresponding dual constraint. 
	By denoting $A_{kl}$ the column of $A$ corresponding to variable $x_{kl}$ and for a feasible dual solution $u$ of $(\dual)$, it is defined as :
	$$ r_{kl}(u) = c_{kl} - uA_{kl}$$
	
	For a pair $kl\in B$, we denote by $R_{kl}$, the exact gap between  $z^*$ (optimal value of $(\primal)$) and $z_{|kl}^*$ (optimal value of $(\primal_{|kl})$) so that $R_{kl} = z_{|kl}^* - z^*$. It is the exact increase of the objective when forcing the pair $kl$ in the solution. A reduced cost $r_{kl}(u)$ of a pair $kl$ is said to be exact and referred to as an \defined{exact reduced cost} when $$r_{kl}(u) = R_{kl} = z_{|kl}^* - z^*$$
	Note that by strong duality, we equivalently have $R_{kl} = w_{|kl}^* - w^*$. By Assumption \ref{feasAssump}, $R_{kl}$ is well defined and finite.
	\begin{example}[Continued: Reduced costs]
		\slshape\sffamily~
		Considering a feasible solution $u$ of the dual formulations $(\dual_{\constraint{WAD}})$ and $(\dual_{\constraint{SP}})$ (more precisely for $(\dual_{\constraint{WAD}})$, $u=(u^1,u^2)$), the reduced costs are given by:
		\begin{multicols}{2}
			For	\constraint{MinWAllDiff} and $(\dual_{\constraint{WAD}})$, 
			$$r_{ij}(u) = c_{ij} -u^1_i -u^2_j$$
			
			\columnbreak
			For \constraint{ShortestPath} and $(\dual_{\constraint{SP}})$, $$\begin{cases}
				r_{ij}(u) =  c_{ij}-u_i+u_j & \forall ij \in \mathcal{A}, j \neq t\\
				r_{it}(u) = c_{it}-u_i & \forall it \in \mathcal{A}
			\end{cases}$$
		\end{multicols}		
	\end{example}
	Reduced costs provide lower bounds of the increase of $z^*$ when variable $x_{ij}$ is forced to one. Since this is a corner stone of the filtering techniques based on LP and the present work, Property \ref{prop: bounds of RC} states it explicitly.
	\begin{property}
		\label{prop: bounds of RC}
		For any optimal solution $u^*$ of $(\dual)$ and any pair $kl \in B$, we have $$0\le r_{kl}(u^*)\le R_{kl}$$
	\end{property}

	\begin{proof}
			The first inequality $r_{kl}(u^*) \geq 0$ is ensured by the corresponding dual constraint since $u^*$ is a feasible solution of $(\dual)$ and the reduced cost is the associated slack. Inequality $r_{kl}(u^*)\le R_{kl}$ is a consequence of sensitivity analysis and we refer the reader to \cite{Wolsey:1998} for more details. A proof is also available in the Annex for an interested reader.
	\end{proof}

	The previous property is the basis for \textit{variable fixing} \citep{Nemhauser:1988:ICO:42805}. 
	Since a reduced cost can be lower than the exact reduced cost, this technique gives an incomplete filtering. Let's give an example for our two illustrative constraints.

	\begin{example}[continued: reduced cost filtering]
	\label{example: rcfiltering}
		\sffamily\slshape
		
		{
		Example of domains (pairs $B$) are represented as graphs below for \constraint{MinWAllDiff} ((a) on the left) and the \constraint{ShortestPath} ((b) on the right). A dual solution is also given with the dual value (in red) associated to each vertex. An edge (for (a)) and an arc (for (b)) represents a pair. Each pair is labelled with its original cost $c_{ij}$ as well as its reduced cost $r_{ij}(u)$ in the proposed dual solutions. Consider $\overline{Z} = 1$ (for both):
		\begin{multicols}{2}
			\begin{tikzpicture}
				\tikzset{
					node/.style={circle, draw, fill=lightgray, minimum size=7mm, inner sep=0pt,scale=.8},
					found/.style={very thick, black!40},
				}
				\def\h{-11mm}
				\def\l{30mm}
				\foreach \x in {0,1,2} {
					\node[node] (x\x) at (0,{\x*\h}) {$X_\x$};
					\node[node] (\x) at (\l,{\x*\h}) {$\x$};
				} 
				\draw[red]
					(x0) node[left=10pt]{$0$} (0) node[right=10pt]{$0$}
					(x1) node[left=10pt]{$1$} (1) node[right=10pt]{$-1$}
					(x2) node[left=10pt]{$3$} (2) node[right=10pt]{$-3$};
					
				\draw ({.5*\l},{-1.5*\h}) node{(a)};
				\draw ({.2*\l},{-\h}) node[node]{$X_i$} node[red,left=10pt] {$u^1_i$} -- node[above=-3pt,scale=.8] {$\left(c_{ij},r_{ij}(u)\right)$} ++({.6*\l},0) node[node] {$j$} node[red,right=10pt] {$u^2_j$} ;
				
				\draw (x0) -- node[scale=.8,above=-3pt]{$(0,0)$} (0);					
				\draw (x0) -- node[pos=.3,scale=.8,above=-3pt,sloped]{$(1,2)$} (1);
				\draw (x1) -- node[pos=.2,scale=.8,above=-3pt,sloped]{$(2,1)$} (0);
				\draw (x1) -- node[pos=.33,scale=.8,above=-3pt,sloped]{$(0,0)$} (1);
				\draw (x1) -- node[pos=.77,scale=.8,above=-3pt,sloped]{$(1,3)$} (2);
				\draw (x2) -- node[pos=.7,scale=.8,above=-3pt,sloped]{$(2,0)$} (1);
				\draw (x2) -- node[scale=.8,above=-3pt,sloped]{$(0,0)$} (2);
			\end{tikzpicture}\par
		
		\columnbreak
			\begin{tikzpicture}
			\tikzset{
			    node/.style={circle, draw, fill=lightgray, minimum size=7mm, inner sep=0pt,scale=.75},
				found/.style={very thick, black!40},
			}
			\def\l{10mm}
			\def\h{10mm}
					
			\draw ({2*\l},{2.5*\h}) node{(b)};
			\draw[node] (0,{2.5*\h}) node[node](i){$i$} node[red,left=7pt] {$u_i$}  ++({2.5*\l},0) node[node](j){$j$} node[red,right=7pt] {$u_j$} ;			
			\draw[-{Latex}] (i) -- node[above=-2.5pt, scale=.8] {$\left(c_{ij},r_{ij}(u)\right)$} (j);
			
			\draw(0,0) node[node](s){$s$} node[left=10pt,red]{0};					
			\draw (\l,\h) node[node](1){1} node[above=7pt, red]{0};
			\draw ({2*\l},0) node[node](2){2}  node[above=7pt, red]{0};
			\draw ({1*\l},-\h) node[node](3){3}  node[below left=5pt, red]{-1};
			\draw ({3*\l},-\h) node[node](4){4}  node[below right=5pt, red]{1};
			\draw ({4*\l},0) node[node](t) {$t$}  node[right=10pt, red]{0};
			\draw[-{Latex}] (s) -- node[pos=.4,above=-2pt,sloped, scale=.8]{$(0,0)$} (1);
			\draw[-{Latex}] (s) -- node[pos=.4,above=-3pt,sloped, scale=.8]{$(2,2)$} (2);
			\draw[-{Latex}] (s) -- node[pos=.4,below=-2pt,sloped, scale=.8]{$(1,0)$} (3);
			\draw[-{Latex}] (1) -- node[above=-3pt,sloped, scale=.8]{$(2,2)$} (t);
			\draw[-{Latex}] (1) -- node[pos=.4,below=-2pt,sloped, scale=.8]{$(0,0)$} (2);
			\draw[-{Latex}] (3) -- node[below=-2pt,sloped, scale=.8]{$(0,2)$} (4);	
			\draw[-{Latex}] (2) -- node[below=-2pt,sloped, scale=.8]{$(0,0)$} (t);	
			\draw[-{Latex}] (4) -- node[pos=.4,below=-2pt,sloped, scale=.8]{$(1,0)$} (t);	
		\end{tikzpicture}\par
		\end{multicols}		
		}	
	 
		\begin{description}
			\item[(a):] $r_{(0,1)}(u)  = 2 > \overline{Z}$ and  $r_{(1,2)}(u) = 3 > \overline{Z}$.			
			Thus, there's no assignment of cost lower than $\overline{Z}$ containing one of the pairs $(0,1)$ or $(1,2)$.
			
			Remark 1: Note that, one of the reduced costs is exact ($r_{(1,2)}(u) = R_{(1,2)}$), whereas the other is not ($r_{(0,1)}(u) < R_{(0,1)} =3$).
			
			Remark 2: $R_{(1,0)} = 3$ thus pair $(1,0)$ is inconsistent event though its reduced cost in $u$ is not high enough to detect it.
			\item[(b):]
			$r_{(s,2)}(u) = 2 > \overline{Z}$ ; $r_{(1,t)}(u) = 2 > \overline{Z}$ and  $r_{(3,4)}(u) = 2 > \overline{Z}$.			
			Thus there's no path of cost lower than $\overline{Z}$ passing through one of the arcs $(s,2)$, $(1,t)$ or $(3,4)$.
			
			Remark: $R_{(s,3)} = 2$ thus the arc $(s,3)$ is inconsistent but $r_{(s,3)}(u) = 0$ so this dual solution doesn't filter this value.
		\end{description}
	\end{example}

Finally, we recall the Complementary Slackness Theorem, which is key to linear programming duality and central to the present work.
We give its statement (from \cite{chvatal1983linear}) using a general linear programming model where the primal and dual problem are respectively denoted $(\primal)$ and $(\dual)$:
$$(\primal) \qquad z^* = \min \{cx : Ax \geq b, x \in \mathbb{R}^n_+\}$$
$$(\dual) \qquad w^* = \max \{ub : uA \leq c, u \in \mathbb{R}^m_+\}$$
Below, index $j$ is used for variables (columns) and $A_j$ is $j$-th column of $A$ whereas index $i$ refers to constraints (rows) and $A^{\intercal}_i$ is the $i$-th line of $A$. 
\begin{theorem}[Complementary Slackness (CS) Theorem] Let $x^*$ be a feasible solution of $(\primal)$ and $u^*$ be a feasible solution of $(\dual)$. Necessary and sufficient conditions for simultaneaous optimality of $x^*$ and $u^*$ are\\
	\indent$u^*A_j = c_j \textrm{ or } x^*_j = 0\qquad (\textrm{or both}) \qquad \forall j \in \{1,\ldots, n\} \qquad  \qquad(\textrm{CS1})$\\
	\indent\indent and\\
	\indent$A^{\intercal}_ix^* = b_i \textrm{ or } u^*_i = 0\qquad (\textrm{or both}) \qquad \forall i \in \{1,\ldots, m\} \qquad  \qquad(\textrm{CS2})$
\end{theorem}
We will use the CS conditions excluding a primal variable $x_k$ and refer to this restricted form of the CS conditions as $k$-excluded. To be precise:
\begin{definition} [$k$-excluded Complementary Slackness (CS) conditions] \label{kCSconditions} Let $x^*$ be a feasible solution of $(\primal)$ and $u^*$ be a feasible solution of $(\dual)$. The CS conditions excluding primal variable $x_k$ are\\
    \indent$u^*A_j = c_j \textrm{ or } x^*_j = 0\qquad (\textrm{or both}) \qquad \forall j \in \{1,\ldots, n\}\textrm{\boldmath$\setminus\{k\}$} \:\:\:\:\: (\textbf{k-CS1})$\\
    \indent\indent and\\
    \indent$A^{\intercal}_ix^* = b_i \textrm{ or } u^*_i = 0\qquad (\textrm{or both}) \qquad \forall i \in \{1,\ldots, m\} \qquad  \qquad(\textrm{CS2})$
\end{definition}

	In the following, we show how to find a set of dual solutions that gives the exact reduced costs, to perform a complete filtering (\ie achieve AC). 
\section{Analysis}
\label{section: analysis}
	For a weighted constraint $W\!C$ and under Assumption \ref{feasAssump}, AC requires identifying each pair $kl$ for which all the integer solutions of $\left(\primal_{|kl}\right)$ have a cost greater than the fixed upper bound $\overline{Z}$. Since $(\primal)$ is assumed ideal, the optimal value $z_{|kl}^*$ of the linear relaxation $\left(\primal_{|kl}\right)$ can be used to establish if value $l$ is consistent for $X_k$ \emph{i.e} if $z^*_{|kl} \leq \overline{Z}$. The assumption that $(\primal)$ is ideal is required to make sure there is no integrality gap and AC can be achieved only by using the linear relaxation $(\primal)$.
	Since the optimal value $z_{|kl}^*$ can be computed as $z^* + R_{kl}$, we are interested in exact reduced costs and whether they can be obtained from dual solutions of $(\primal)$. 
	
	
	Two pairs $ij$ and $kl$ are said \defined{incompatible} if they do not belong together to a support. More precisely, there is no support for $C$ with $X_i = j$ and $X_k = l$. The set $\set{I} \subseteq B$ refers to a set of pairs that are (simultaneously) pairwise incompatible. We start by defining a modified problem denoted $\left(\tilde{\primal}_{\set{I}}\right)$ that is repeatedly used throughout the document. It is instrumental as we can show that there exists an optimal solution of its dual, $\left(\tilde{\dual}_{\set{I}}\right)$, with exact reduced costs of all pairs of $\set{I}$. In particular, it shows the first important step (Property \ref{prop: exists u st r_kl=R_kl}) that for any pair $kl$ there exists an optimal solution of $(\dual)$ with the exact reduced cost of $kl$. This property justifies that a complete filtering (\ie achieve AC) is possible using reduced costs.
	
	\begin{definition}[Modified problem $\left(\tilde{\primal}_{\set{I}}\right)$]
	\label{defPtilde}
	Let $\set{I}$ be a set of pairwise incompatible pairs of $B$ and let Assumption \ref{feasAssump} hold. Problem $\left(\tilde{\primal}_{\set{I}}\right)$ is the formulation identical to $\left(\primal\right)$ except for the costs related to the pairs of \set{I}. More precisely,
		$$\tilde{c}_{ij} = \begin{cases}
			c_{ij}-R_{ij}\quad &\forall ij\in \set{I}\\
			 c_{ij} &\forall ij\in B\setminus \set{I}
		\end{cases}$$
	\end{definition}
	
	\begin{remark}[]
	\label{feasAndBoundedPtilde}
Under assumption \ref{feasAssump}, formulation $\left(\tilde{\primal}_{\set{I}}\right)$ is feasible, bounded and ideal.
	\end{remark}
	\begin{proof}
	 $\left(\tilde{\primal}_{\set{I}}\right)$ is feasible since $\left(\primal\right)$ is feasible. It is bounded since all $R_{ij}$ are finite for any $ij$ of $\set{I}$ and the domains are finite. Since only the costs are changed, the extreme points of $\left(\tilde{\primal}_{\set{I}}\right)$ are unchanged.
	\end{proof}
	\begin{property}[]
	\label{optSolDtilde}
    For any set of pairwise incompatible pairs $\set{I}$, there exists an optimal dual solution $u^*$ of $\left( \tilde{\dual}_{\set{I}}\right)$ that is optimal for $\left(\dual\right)$ 
    and such that all reduced costs of the pairs of $\set{I}$ are exact: $$r_{ij}(u^*) = R_{ij} \:\:\forall ij \in \set{I}$$
	\end{property}
	\begin{proof}
		Let $\tilde{u}^*_{\set{I}}$ be an optimal solution for $\left(\tilde{\dual}_{\set{I}}\right)$, the dual of $\left(\tilde{\primal}_{\set{I}}\right)$, and $\tilde{z}^*_{\set{I}}$ its value. Such a solution exists by Observation \ref{feasAndBoundedPtilde}. We show that $\tilde{u}^*_{\set{I}}$ is feasible and optimal for $\left(\dual\right)$ while providing the exact reduced costs of all pairs of $\set{I}$:
		\begin{itemize}
			\item For any pair $ij$ of $\set{I}$, since $R_{ij}\ge 0$ and $\tilde{c}_{ij}\le c_{ij}$ we have:\\
			$\begin{cases}
				\tilde{u}^*_{\set{I}}A_{ij} \le \tilde{c}_{ij}\le c_{ij} \quad \forall ij \in \set{I}\\
				\tilde{u}^*_{\set{I}}A_{ij}\le \tilde{c}_{ij}= c_{ij}\quad \forall ij \in B\setminus \set{I}
			\end{cases}
			$\\
			and $\tilde{u}^*_{\set{I}}$ is a feasible solution for $\left(\dual\right)$.
			\item 
			Since costs are lower in $\left(\tilde{\primal}_{\set{I}}\right)$, we know that $\tilde{z}^*_{\set{I}} \leq z^*$. Let's assume that $\tilde{z}^*_{\set{I}} < z^*$. At least one pair of $\set{I}$ must be used in such an optimal solution $\tilde{x}^*$ of $\left(\tilde{\primal}_{\set{I}}\right)$ otherwise we would have $\tilde{z}^*_{\set{I}} = z^*$. Because $\left(\tilde{\primal}_{\set{I}}\right)$ is an ideal formulation, there is an integer optimal solution and we can assume $\tilde{x}^*$ is integral. Additionally, since $\set{I}$ is a set of incompatible pairs,  all $\tilde{x}^*_{ij}$ for $ij\in \set{I}$ are equal to zero except one, $\tilde{x}^*_{kl}=1$. The objective value of $\tilde{x}^*$ in  $\left(\primal\right)$ is therefore $\tilde{z}^*_{\set{I}} + R_{kl} 
				= \tilde{z}^*_{\set{I}} + z_{|kl}^* - z^* < z^* + z_{|kl}^* - z^* = z_{|kl}^*$. That is impossible by definition of $z_{|kl}^*$ so  $\tilde{z}^*_{\set{I}} =z^*$ and $\tilde{u}^*_{\set{I}}$ is optimal for $\left(\dual\right)$.
			\item Since $\tilde{u}^*_{\set{I}}$ is feasible for $\left(\tilde{\dual}_{\set{I}}\right)$, for any $ij \in \set{I}$, we have $\tilde{u}^*_{\set{I}}A_{ij} \leq \tilde{c}_{ij}$ \emph{i.e}
			$\begin{aligned}[t]
				&& \tilde{u}^*_{\set{I}}A_{ij}& \leq c_{ij} - R_{ij}\\
				\iff && c_{ij} - \tilde{u}^*_{\set{I}}A_{ij} &\ge R_{ij}\\
				\iff && r_{ij}(\tilde{u}^*_{\set{I}}) &\ge R_{ij}\\
			\end{aligned}$\\
			Because $\tilde{u}^*_{\set{I}}$ is also an optimal solution for $\left(\dual\right)$, from Property \ref{prop: bounds of RC}, ${r_{ij}(\tilde{u}^*_{\set{I}}) \leq R_{ij}}$ and therefore $r_{ij}(\tilde{u}^*_{\set{I}}) = R_{ij} \:\: \forall ij \in \set{I}$.
	    \end{itemize}
    \end{proof}
	
	{\begin{property}[Existence of an optimal dual solution $u^*$ \st ${r_{kl}(u^*) = R_{kl}}$]
		\label{prop: exists u st r_kl=R_kl}		
		Under assumption \ref{feasAssump}, for any pair $kl\in B$
		, there exists an optimal dual solution $u^*$ of $(\dual)$ such that $r_{kl}(u^*) = R_{kl}$.
	\end{property}
	\begin{proof}
	The solution $u^*$ can be built by Property \ref{optSolDtilde} as an optimal solution of $\left(\tilde{\dual}_{\set{I}}\right)$ 
	with $\set{I}$ restricted to the single pair $kl$: $\set{I} = \{kl\}$. 
	\end{proof}	

	In practice, as we will see below, a single dual solution might provide the exact reduced costs of many pairs.
	
	\medskip
	The result of \citep{german2017arc} is based on an interior point of the primal. We now give a similar result in the dual by considering that a dual solution with strictly positive reduced costs (slacks) can be seen as an "interior point of the dual". A constraint $C$ can be encoded as a weighted constraint $W\!C$ with 0/1 costs. Let $Z$ be equal to 0. We assume that the domains of the variables in $W\!C$ can be extended by adding values so that Assumption \ref{feasAssump} is satisfied. For instance, by considering a complete bipartite graph for $\textsc{AllDifferent}$. An original value $j$ of a domain $D(X_i)$ is given a cost of 0 ($c^{01}_{ij} = 0$) and the remaining pairs $ij$ (encoding values \textbf{not} present in the initial domains) are given a cost of 1 ($c^{01}_{ij} = 1$).
	A pair belongs to a support of $C$ if and only if it belongs to a support of cost 0 in $WC$. 
	Property \ref{prop: r_ij>0 <=> ij filterable} shows that any positive reduced cost exhibits an inconsistent pair and a single dual solution can rule out all inconsistent pairs (values).
	Let $\set{F}$ be the set of inconsistent pairs \ie the pairs that do not belong to any support for $C$.
	
	\begin{property}[AC for $C$ with a single dual solution]
		\label{prop: r_ij>0 <=> ij filterable}
		Let $W\!C$ be the 0/1 encoding of $C$ under Assumption \ref{feasAssump}.
		There exists an optimal dual solution $\tilde{u}$ of $(\dual)$ for $W\!C$ with costs $c^{01}$ such that for any pair $ij$ of $B$:
		 $$ij \in \set{F} \iff r_{ij}(\tilde{u}) >0$$ 
	\end{property}
\begin{proof}
		 
		The cost encoding $c^{01}$ implies $R_{ij}\ge 1$ $\forall ij\in\set{F}$. We can consider a set of optimal solutions of $(\dual)$:  $\left\{\tilde{u}^{ij}: ij \in \set{F}, r_{ij}(\tilde{u}^{ij}) \ge 1\right\}$. Property \ref{prop: exists u st r_kl=R_kl} ensures that these solutions exist. One can remark that $\forall kl \not \in \set{F}$, $R_{kl} = 0$ and thus, $\forall ij\in \set{F},~r_{kl}(\tilde{u}^{ij}) = 0$. Let $\tilde{u}$ be the average solution of the previous set: $$\tilde{u} = \frac{1}{\left|\set{F}\right|} \sum\limits_{ij\in\set{F}}\tilde{u}^{ij}$$
		This solution is feasible, optimal, and $r_{ij}(\tilde{u}) > 0$ for all $ij \in \set{F}$ and  $r_{ij}(\tilde{u}) = 0$ otherwise ($\forall ij \not \in \set{F}$).
	\end{proof}


	Let's go back to the weighted case $W\!C$.
	A single dual solution can, in fact, exhibit the exact reduced costs of many pairs.  Considering a given variable $x_{kl}$ of $\left(\primal\right)$, the $kl$-excluded CS conditions (see Definition \ref{kCSconditions}), are necessary and sufficient conditions for a feasible solution of $\left(\primal_{|kl}\right)$ and an optimal solution of $\left(\dual \right)$ to respectively be optimal for $\left(\primal_{|kl}\right)$ and give the exact reduced cost of $kl$. This result is analogous to the CS theorem but characterizes optimality regarding a given pair. 
	\begin{theorem}[AC complementary slackness]
		\label{ac slackness}
		 Let $x^*$ be a feasible solution of $\left(\primal_{|kl}\right)$ and $u^*$ an optimal dual solution of $\left(\dual\right)$. The solution $x^*$ is optimal for $\left(\primal_{|kl}\right)$ and $r_{kl}(u^*) = R_{kl}$ if and only if $x^*$ and $u^*$ satisfy the $kl$-excluded complementary slackness conditions.
	\end{theorem}
    \begin{proof} Note that the dual $\left(\dual_{|kl}\right)$ has an additional dual variable $\alpha_{kl}$ related to the constraint $x_{kl}=1$. Since $u^*$ is a feasible solution of $(\dual)$, a dual feasible solution for $\left(\dual_{|kl}\right)$ is readily available by using the values of $u^*$ and setting $\alpha^*_{kl}$ to the slack $\alpha^*_{kl}=r_{kl}(u^*)$.
    
    $\Leftarrow$: Let's show that the dual solution $(u^*,\alpha^*_{kl}=r_{kl}(u^*))$ is optimal for $\left(\dual_{|kl}\right)$ and that $x^*$ is optimal for $\left(\primal_{|kl}\right)$ by checking they satisfy the CS conditions. Firstly, the CS conditions are satisfied for all variables except $x_{kl}$ by assumption ($kl$-excluded CS are assumed). Secondly, (CS1) is satisfied for $x_{kl}$ by construction since we fill the slack with $\alpha^*_{kl}$. Finally, the additional constraint $x_{kl}=1$ is tight so (CS2) remains valid. Thus, according to the complementary slackness theorem, $x^*$ is optimal for $\left(\primal_{|kl}\right)$.
    Its cost is therefore $z^*_{|kl}$. The corresponding optimal dual solution $(u^*,\alpha^*_{kl}=r_{kl}(u^*))$ of $\left(\dual_{|kl}\right)$ has a cost of $w^*_{|kl} = w^* + \alpha^*_{kl}= w^* + r_{kl}(u^*)$. By strong duality ($w^*_{|kl} = z^*_{|kl}$ and $w^* = z^*$), we have $r_{kl}(u^*) = z^*_{|kl} - z^*= R_{kl}$.
    
	$\Rightarrow$: Assuming $r_{kl}(u^*) = R_{kl}$, the dual solution $(u^*_{kl},\alpha^*_{kl}=r_{kl}(u^*))$ is optimal for $\left(\dual_{|kl}\right)$ because it is feasible for $\left(\dual_{|kl}\right)$ and its objective is $w^* + R_{kl} = z^*_{|kl}$. Since $x^*$ is assumed to be optimal for $\left(\primal_{|kl}\right)$, it satisfies the CS conditions with $(u^*,\alpha^*_{kl}=r_{kl}(u^*))$. As a result, $x^*$ and $u^*$ satisfy $kl$-excluded CS conditions since (CS2) are unchanged.
	\end{proof}

	This theorem shows how to use complementary slackness to know which reduced costs are exact given an optimal dual solution. Typically, if there exists a primal solution (a support) using $kl$ (with $x_{kl} = 1$) and satisfying $kl$-excluded conditions with $u^*$ then the reduced cost for $kl$ is exact in $u^*$. Note that when the primal problem is made of equality constraints, conditions (CS2) of the complementary slackness is true for any feasible solution since there is no possible slack. This is the case for our two examples related to assignments and paths. As a result, we know that $r_{kl}(u^*)$ is exact if $kl$ belongs to a support of $C$ where the other pairs (all pairs except $kl$) have a null reduced costs in $u^*$. 
	Characterizing exact reduced costs boils down to a feasibility problem using the values of null reduced costs.
	
	\begin{example}[Characterization of exact reduced costs]
		\label{ex: characterization r_kl=R_kl}
		\slshape\sffamily~
		
		Consider $\overline{Z} = 1$. We apply the previous theorem by identifying the exact reduced costs of a given dual optimal solution of $(\dual)$. The two cases of \constraint{MinWAllDiff} and \constraint{ShortestPath} are illustrated. Note that the two instances differ from the ones of Example \ref{example: rcfiltering}.
		\begin{multicols}{2}
			\begin{tikzpicture}
				\tikzset{
					every node/.style={scale=.8},
					node/.style={circle, draw, fill=lightgray, minimum size=7mm, inner sep=0pt,scale=1.2},
					support/.style={very thick},
				}
				\def\h{-11mm}
				\def\l{30mm}
				
				\draw ({.5*\l},{-1.5*\h}) node{(a)};
				\foreach \x in {0,1,2} {
					\node[node] (x\x) at (0,{\x*\h}) {$X_\x$};
					\node[node] (\x) at (\l,{\x*\h}) {$\x$};
				} 
				\draw[red]
				(x0) node[left=10pt]{$0$} (0) node[right=10pt]{$0$}
				(x1) node[left=10pt]{$-1$} (1) node[right=10pt]{$1$}
				(x2) node[left=10pt]{$0$} (2) node[right=10pt]{$0$};
				
				\draw ({.2*\l},{-\h}) node[node]{$X_i$} node[red,left=10pt] {$u^1_i$} -- node[above=-3pt] {$\left(c_{ij},r_{ij}(u)\right)$} ++({.6*\l},0) node[node] {$j$} node[red,right=10pt] {$u^2_j$} ;
				
				\draw (x0) -- node[above=-3pt]{$(0,0)$} (0);					
				\draw (x0) -- node[pos=.32,above=-4pt,sloped]{$(2,1)$} (1);
				\draw[support] (x0) -- node[pos=.12,below=-3pt,sloped]{$(0,\boldsymbol{0})$} (2);
				\draw[support,dash pattern={on 3pt off 1pt}] (x1) -- node[pos=.7,above=-4pt,sloped,black]{$(1,\boldsymbol{2})$} (0);
				\draw (x1) -- node[pos=.65,above=-4pt,sloped]{$(0,0)$} (1);
				\draw (x1) -- node[pos=.7,below=-3pt,sloped]{$(1,2)$} (2);
				\draw[support] (x2) -- node[pos=.15,above=-4pt,sloped]{$(1,\boldsymbol{0})$} (1);
				\draw (x2) -- node[below=-3pt,sloped]{$(0,0)$} (2);
			\end{tikzpicture}\par
		
		\columnbreak
		
			\begin{tikzpicture}
				\tikzset{
					every node/.style={scale=.8},
					node/.style={circle, draw, fill=lightgray, minimum size=7mm, inner sep=0pt},
					support/.style={very thick},
				}
				\def\l{10mm}
				\def\h{10mm}
				
				\draw ({2*\l},{2.5*\h}) node{(b)};
				\draw(0,0) node[node](s){$s$} node[left=10pt,red]{0};					
				\draw (\l,{1.2*\h}) node[node](1){1} node[above=10pt, red]{$-1$};
				\draw ({2*\l},0) node[node](2){2}  node[above right=7pt, red]{$-1$};
				\draw ({2*\l},-\h) node[node](3){3}  node[below=10pt, red]{0};
				\draw ({4*\l},0) node[node](t) {$t$}  node[right=10pt, red]{0};
				
				\draw[-{Latex},support] (s) -- node[pos=.4,above=-2pt,sloped]{$(1,0)$} (1);
				\draw[-{Latex}] (s) -- node[pos=.4,above=-2pt,sloped]{$(2,1)$} (2);
				\draw[-{Latex}] (s) -- node[pos=.4,below=-2pt,sloped]{$(0,0)$} (3);
				\draw[-{Latex}] (1) -- node[above=-2pt,sloped]{$(1,2)$} (t);
				\draw[-{Latex},support] (1) -- node[pos=.4,below=-2pt,sloped]{$(0,0)$} (2);
				\draw[-{Latex},support] (2) -- node[pos=.4,below=-2pt,sloped]{$(1,2)$} (t);	
				\draw[-{Latex}] (3) -- node[below=-2pt,sloped]{$(0,0)$} (t);	 (t);			
				
				\draw[node] (0,{2*\h}) node[node](i){$i$} node[red,left=10pt] {$u_i$}  ++({2*\l},0) node[node](j){$j$} node[red,right=10pt] {$u_j$} ;			
				\draw[-{Latex}] (i) -- node[above=-3pt] {$\left(c_{ij},r_{ij}(u)\right)$} (j);
			\end{tikzpicture}\par
		\end{multicols}
		\begin{description}
			\item[(a):]
			$\set{S} =\left\{(0,0), (1,1), (2,2) \right\}$ is a support of minimal cost (its cost is 0).
			
			 $\set{S}' =\left\{(0,2), (1,0) , (2,1) \right\}$ is a support (of minimal cost) containing $(1,0)$ \ie a perfect matching of minimum weight using pair $(X_1,0)$. Moreover the reduced costs of two pairs out of three are null: $r_{(0,2)}(u^*) = r_{(2,1)}(u^*) = 0$. 			
			Thus $r_{(1,0)}(u^*)$ is exact \ie $R_{(1,0)} = r_{(1,0)}(u^*) = 2$.
			\item[(b):] 
			$\set{S} =\left\{(s,1) ; (1,2) ;(2,t) \right\}$ is a support (of minimal cost) containing $(2,t)$ with $r_{(s,1)}(u^*) = r_{(1,2)}(u^*) = 0$. 			
			Thus $R_{(2,t)} = r_{(2,t)}(u^*)=2$. In other words, a shortest path using pair $(2,t)$ has a cost of $2$.
		\end{description}
	\end{example}

	For the specific \constraint{MinWAllDiff} constraint, the conditions of the previous theorem are met if and only if there exists a cycle, alternating with respect to a support of minimal cost, such that all its reduced costs are 0 except one which is exact. Similarly for \constraint{ShortestPath}, if a path from $s$ to $t$ can be built using the arcs of the null reduced costs except for a single additional arc, then the reduced cost of this additional arc is exact. 
	
	Note that the \textbf{optimality} of the reduced cost is checked by solving a \textbf{feasibility} problem. This latter problem is stated by distinguishing the pairs of null reduced costs from the remaining pairs. In particular, it does not use the precise value of the costs themselves. Very similarly, optimality is reached in a primal dual algorithm when a feasible solution is obtained with the pairs of null reduced costs alone. For instance, the Hungarian algorithm stops when a maximum matching of the graph of null reduced costs has a cardinality of $n$. The costs are not used when checking this condition, they have been \emph{combinatorialized} as explained in \citep{papadimitriou1998combinatorial}. From this point of view, we believe Theorem \ref{ac slackness} extends complementary slackness very naturally to deal with arc-consistency. The statement below is a direct consequence of Theorem \ref{ac slackness} but explicitly states necessary and sufficient conditions for a reduced cost to be exact which makes it easier to manipulate.	

		\begin{corollary}[characterization of $r_{kl}(u^*)=R_{kl}$]
	\label{prop: characterization r_kl=R_kl}
		
		Let $u^*$ be an optimal dual solution of $\left(\dual\right)$ and $kl$ a pair of $B$. The reduced cost $r_{kl}(u^*)$ is exact (${r_{kl}(u^*) = R_{kl}}$) if and only if there exists a feasible solution of $\left(\primal_{|kl}\right)$ that satisfies $kl$-excluded complementary slackness conditions with $u^*$.

	\end{corollary}
	
    \begin{proof}
		$\Rightarrow$: we can build explicitly the solution $x^*$ with $x^*_{kl} = 1$ that satisfies $kl$-excluded CS conditions with $u^*$. It can be found as an optimal solution $x^*$ of the modified problem $\left(\tilde{\primal}_{\set{I}}\right)$ with $\set{I} = \{kl\}$ (Definition \ref{defPtilde}). Such a solution exists by Observation \ref{feasAndBoundedPtilde} and we know, as a consequence of Property \ref{optSolDtilde} that it is optimal for $\left(\primal\right)$ so $\tilde{z}^*_{\set{I}} = z^*$. Note that solution $u^*$ is also feasible and optimal for $\left(\tilde{\dual}_{\set{I}}\right)$. The reason is that $\tilde{z}^*_{\set{I}} = z^*$ and $u^*$ is feasible for $\left(\tilde{\dual}_{\set{I}}\right)$ because we assumed $r_{kl}(u^*) = R_{kl}$. Thus $x^*$ and $u^*$ satisfy CS conditions as optimal solutions of $\left(\tilde{\primal}_{\set{I}}\right)$ and $\left(\tilde{\dual}_{\set{I}}\right)$. As a consequence, they satisfy $kl$-excluded conditions as (respectively) feasible solution of ($\primal$) and optimal of ($\dual$).\\
		$\Leftarrow$: By Theorem \ref{ac slackness}, if we have a feasible solution $x^*$ with $x^*_{kl} = 1$ satisfying $kl$-excluded CS with $u^*$ then $r_{kl}(u^*) = R_{kl}$.
	\end{proof}
	We can now turn our attention to dual solutions providing multiple exact reduced costs at the same time.
	For a given optimal solution $u$ of $(\dual)$, let's denote by $\set{R}_u=\left\{ij\in B ~|~r_{ij}(u)=R_{ij}\right\}$, the set of pairs whose reduced costs are exact with respect to $u$. Consider a family $\mathcal{T}$ of optimal solutions of $(\dual)$ denoted $\left\{u^t\right\}_{t\in \set{T}}$. We denote by $E(\mathcal{T}) = \bigcup_{t\in \set{T}} \set{R}_{u^t}$, the set of pairs whose exact reduced costs are provided by $\set{T}$. 
	\begin{definition}[Complete family of dual solutions]
	\label{completeSetDEF}
	A family $\set{T}$ of optimal dual solutions $\left\{u^t\right\}_{t\in \set{T}}$ of $(\dual)$ is said to be \textbf{complete with respect to a set $H$} of pairs if $H \subseteq E(\mathcal{T})$. 
	
	A family is said \textbf{complete} when it is complete with respect to $B$.
	\end{definition}
	In a number of cases, where some dependencies among values hold, AC can be achieved by only computing the exact reduced costs of a subset of pairs $H\subseteq B$. In such a case, AC is derived from a family that is only complete with respect to $H$. This is the case for the example \ref{ex: sets of incompatible edges} below and the path constraint but also for the weighted stable set problem addressed in Section \ref{section: MWIS}.
	To minimize the number of calls to the simplex algorithm in order to compute arc-consistency, we are interested in complete families of minimal cardinality. 
	
	Theorem \ref{ac slackness} implies necessary 
	conditions for two reduced costs to be exact in the same dual solution. Two pairs $kl$ and $ab$ are said to be \textbf{opt-disjoint} if the two sets of optimal solutions of ($\primal_{|kl}$) and ($\primal_{|ab}$) are disjoint. In other words, there is no support containing both $X_k=l$ and $X_a=b$ that would be a minimal support for each of them. This comes from the excluded CS conditions which require all reduced costs to be null except one.
	
	\begin{theorem}[simultaneous exact reduced costs]
		\label{simultaneous exact rc}
		 Let $\set{I}$ be a set of pairs of $B$ with strictly positive exact reduced costs ($R_{kl} > 0 \:\:\forall kl \in \set{I}$). If there exists an optimal solution $u^*$ of $\left(\dual\right)$ with $r_{kl}(u^*) = R_{kl} \:\:\forall kl \in \set{I}$ then pairs of $\set{I}$ are pairwise opt-disjoint.
	\end{theorem}
    \begin{proof} 
		Let's assume that $u^*$ is an optimal solution of $(\dual)$ with exact reduced costs for the pairs of $\set{I}$. By Theorem \ref{ac slackness}, any optimal solution $x^*$ of $\left(\primal_{|kl}\right)$ for a pair $kl \in \set{I}$ must satisfy $kl$-excluded CS conditions with $u^*$. So the reduced costs of all other pairs must be null and such a solution $x^*$ can not include another pair of $\set{I}$ (since for all pairs of $\set{I}$, we have $R_{kl} > 0$). All pairs of $\set{I}$ are therefore pairwise opt-disjoint.
	\end{proof}
	
     The opt-disjoint criterion is necessary but not well suited to derive an algorithm. We already mentioned (Property \ref{optSolDtilde}) a simple sufficient condition that is highlighted below as Corollary \ref{prop: incompatible edges}. Note that two incompatible pairs are opt-disjoint. But incompatibility is sufficient and the resulting condition can be more easily used by an algorithm:
	\begin{corollary}[Incompatible pairs]
		\label{prop: incompatible edges}
		
		For any set $\set{I} \subseteq B$ of pairwise incompatible pairs, there exists an optimal solution of $\left(\dual\right)$ that provides the exact reduced cost of each pair of $\set{I}$.
		
	\end{corollary}
\begin{proof}
		\prooflabel{incompatible edges}
		The dual solution of $(\dual)$ can be built using $(\tilde{D}_{\set{I}})$ (Property \ref{optSolDtilde}).
	\end{proof}

	Domains are simple examples of incompatible pairs since a variable takes a single value. 
	Therefore no two values of its domain can be used together and for each variable $X_k$, the set $\set{F}_k=\left\{kj \in B\right\} = \{kj~|~j \in D(X_k)\}$ is a set of incompatible pairs. Let $u^\star_k$ be an optimal dual solution for which all reduced costs of $\set{F}_k$ are exact. Then, $\set{F} = \left\{u^\star_k ~|~ k \in U\right\}$  is a complete family of dual solutions (recall $U = \{X_1,\ldots, X_n\}$). It is simply based on the domains. There might exist a complete set of dual solutions based on incompatibility that have a smaller cardinality and we illustrate it below with the case of \constraint{ShortestPath}:

%
%
	\begin{example}[Sets of incompatible pairs]
		\label{ex: sets of incompatible edges}
		
		Recall for \constraint{ShortestPath} that variables $x_{ii}$ representing pairs $ii$ (dummy value for $X_i = i$) have not been included in its formulation. We consider a complete family of dual solutions with respect to set $\{ij | ij \in B, i\neq j\} \subset B$. 
		\slshape\sffamily~
		\begin{multicols}{2}
			\begin{minipage}{.9\linewidth}
			For \constraint{MinWeightedAllDiff}, the sets $\set{F}_k$ are sets of incompatible pairs:
			
			\includegraphics[width=4.5cm]{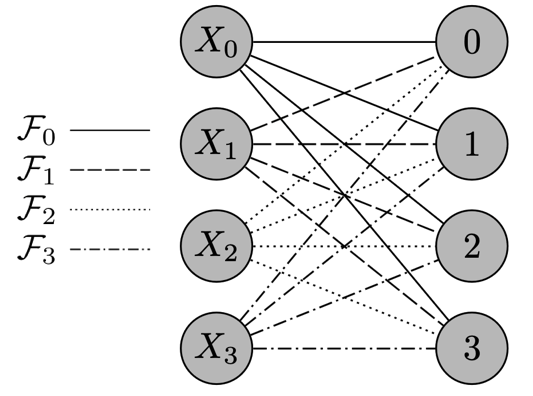}
			
			For each $\set{F}_k$ there exists an optimal dual solution $u^*_k$ of $(\dual_{\constraint{WAD}})$ that gives the exact reduced costs of all pairs in $\set{F}_k$. Thus $\set{F} = \{u^*_k~|~ 0\le k \le 3\}$ is a complete family with 4 dual solutions.
		\end{minipage}
		\columnbreak
		
		\begin{minipage}{.9\linewidth}
		
			For \constraint{ShortestPath}, the topological layers $\set{T}_i$ are sets of incompatible pairs:
			
		    \includegraphics[width=4.5cm]{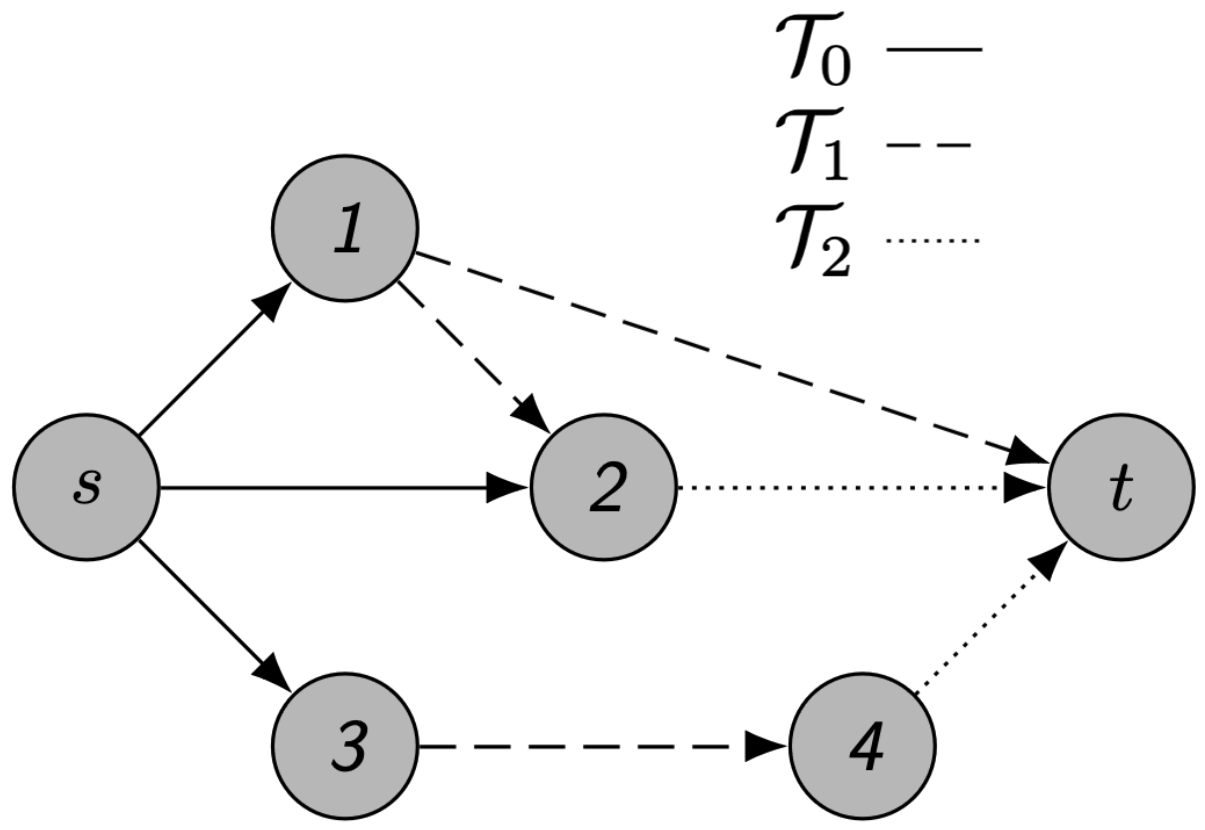}
				
				
			For each $\set{T}_i$ there exists an optimal dual solution $u^*_i$ of $(\dual_{\constraint{SP}})$ that gives all the exact reduced costs for arcs in $\set{T}_i$. Thus $\{u^*_i~|~ 0\le i \le 2\}$ is a complete family of dual solutions of cardinality 3. Note that the complete family $\set{S}$ based on the variables (the vertices) would require 5 dual solutions. 
		\end{minipage}
		\end{multicols}
	
	\end{example}

    An upper bound and a worst case lower bound of this cardinality are given in the following properties:

	\begin{property}[Complete family of $n$ dual solutions]
		\label{prop: complete set of card n}
		
		There exists a complete family of at most $n$ 
		optimal solutions of $(\dual)$.
	\end{property}
	\begin{proof}
		For each $k\in U$, we set $\set{I}_k=\left\{kj \in B\right\}$. Because of the domains, two pairs of $\set{I}_k$ can not belong to the same solution. By Corollary \ref{prop: incompatible edges} there exists an optimal dual solution $\tilde{u}_{\set{I}_k}^*$ that provides all the exact reduced costs of $\set{I}_k$. Since $B = \bigcup\limits_{k\in U} \set{I}_k$, family $\left\{\tilde{u}^*_{\set{I}_k}\right\}_{k\in U}$ of $n$ dual solutions is complete. Note that the solutions are not necessarily distinct.
	\end{proof}

	 This set is possibly not minimal as seen in example \ref{ex: sets of incompatible edges} for \constraint{ShortestPath} where $n=5$ and the set of topological layers gives a complete family of cardinality 3. 
	 Finally, we can show that, in the worst case, the cardinality of such a family is at least $n$. We use the \constraint{MinWAllDiff} to do the proof which was presented in \citep{claus2020wad} but can now be simplified based on Theorem \ref{ac slackness}. As outlined previously, this theorem tells us that a single dual solution can not provide the exact reduced costs of two pairs if they belong together in a common support that is minimum for each. 
	 The result below shows an instance of the assignment problem where a single support is minimum for $n$ values so that $n$ distinct dual solutions are required.
	
	\begin{property}
	\label{prop:upperBoundNbDSol}
		In the worst case, for $W\!C$, $n$ optimal solutions of $(\dual)$ are needed to obtain all exact reduced costs.
	\end{property}
	
	\begin{proof}		
		We consider the \constraint{MinWAllDiff} constraint, and the instance where $$c_{ij} = \begin{cases}
			0 & \text{if $i\ge j$}\\
			1 & \text{otherwise.}
		\end{cases}$$
		
		{
			\centering
			\begin{tikzpicture}
				\tikzset{
					node/.style={circle, draw, fill=lightgray, minimum size=10mm, inner sep=0pt,scale=.7},
				}
				\def\l{30mm}
				\def\h{-10mm}
				\draw
				(0,0) node[node] (x0) {$X_0$} 
				(\l,0) node[node](0){0}
				(0,{\h}) node[node](x1){$X_1$} 
				(\l,{\h}) node[node](1){1} 
				(0,{1.65*\h}) node{$\vdots$}
				(\l,{1.65*\h}) node{$\vdots$}
				(0,{2.5*\h}) node[node](xn-1){$X_{n-1}$} 
				(\l,{2.5*\h}) node[node](n-1){$n-1$} 
				(0,{3.5*\h}) node[node](xn){$X_n$} 
				(\l,{3.5*\h}) node[node](n){$n$}; 	
				\draw 					
				({1.5*\l},-10mm) --  ++(7mm,0) 	({1.5*\l},-10mm) node[above right]{$c_{ij} = 1$ (continuous)}
				(x0) -- (1)
				(x0) --	(n-1)
				(x0) -- (n)
				(x1) -- (n-1)
				(x1) -- (n)
				(xn-1) -- (n);
				\draw[very thick]
				({1.5*\l},-10mm) ++(7mm,0) --  ++(7mm,0)
				({1.5*\l},-20mm)-- ++(7mm,0) ({1.5*\l},-20mm) node[above right]{$\set{L}$ (bold)}
				(x0) -- (n);
				
				\draw[dashed]
				({1.5*\l},0) -- ++(7mm,0) ({1.5*\l},0) node[above right]{$c_{ij} = 0$ (dashed)} 
				(x0) -- (0)
				(x1) -- (0)
				(x1) -- (1)					
				(xn-1) -- (0)
				(xn-1) -- (1)
				(xn-1) -- (n-1)
				(xn) -- (0)
				(xn) -- (1)
				(xn) -- (n-1)
				(xn) -- (n);
				\draw[very thick, dashed] 
				({1.5*\l},0) ++(7mm,0) -- ++(7mm,0) 
				({1.5*\l},-20mm) ++(7mm,0) -- ++(7mm,0) 
				(x1) -- (0)
				(xn-1) -- ++({.3*\l}, {-.3*\h}) 
				(1) -- ++({-.3*\l}, {.3*\h}) 
				(xn) -- (n-1);
				\draw[very thick, dash pattern={on 1pt off 2pt}] 
				(xn-1) ++({.3*\l}, {-.3*\h}) --  ++({.15*\l}, {-.15*\h}) 
				(1) ++({-.3*\l}, {.3*\h}) -- ++({-.15*\l}, {.15*\h}) ;
				
			\end{tikzpicture}\par
		}
		
		Let $\set{L} = \left\{(i,i-1) ~|~ \forall\, 1 \le i\le n\right\}\cup \{(0,n)\}$. Let us show that the reduced costs of the pairs of $\set{L}$ can not be pairwise exact for the same dual solution.
		
		For any pair $kl \in \set{L}$, $R_{kl}=1$ and $\set{L}$ is a support of minimal cost (its cost is 1) containing this pair. From Theorem \ref{ac slackness} if $u^*$ is an optimal dual solution of $(\dual)$ with $r_{kl}(u^*) = R_{kl}$ and $\set{L}$ is an optimal solution of $(\primal_{|kl})$, then it satisfies $kl$-excluded CS conditions so the reduced costs for the pairs of $\set{L}\setminus \{kl\}$ must be null. Thus only one exact reduced cost of \set{L} can be given by an optimal dual solution and at least $\left|\set{L}\right| = n$ dual solutions are needed to obtain all exact reduced costs.
	\end{proof}
	
	For a set of incompatible pairs \set{I}, if we modify the objective function of $(\dual)$ to maximize $\sum_{ij\in \set{I}} r_{kl}(u^*)$ by keeping $u^*$ optimal, we will obtain a solution that provides the exact reduced costs for all pairs of \set{I}. This technique is the foundation of the algorithm proposed in the next section.

\section{An LP based algorithm}
\label{section: LP base algo}

	
	For a given set of incompatibles pairs, $\set{I}$, we can modify $\left(\dual\right)$ to obtain  a solution which is appropriate for filtering.		
	
	$$\left(\dual_{\set{I}}^z\right) \qquad  \max_{u} \{ w^{z}_{\mathcal{I}} = \sum_{kl \in \mathcal{I}} r_{kl}(u) : uA \leq c, ub = z^*, u \in \mathbb{R}^m_+\} $$
		
		All the constraints of $\left(\dual\right)$ are included in $\left(\dual_{\set{I}}^z\right)$ and the additional constraint\, $ub = z^*$ ensures the obtained solution is optimal for $\left(\dual\right)$. 
		Corollary \ref{prop: incompatible edges} shows that there exists an optimal solution of $\left(\dual\right)$ in which all exact reduced costs of $\set{I}$ are reached. Thanks to the property \ref{prop: exists u st r_kl=R_kl}, the objective function, which is the sum of the reduced costs of $\set{I}$, ensures that such a solution is found.
	
	A drawback of $\left(\dual_\set{I}^z\right)$ is the preliminary computation of $z^*$. Constraint ${ub = z^*}$ also considerably changes the formulation of the original dual $\left(\dual\right)$ which might be inconvenient when a dedicated algorithm is available for solving $\left(\dual\right)$. But $\left(\dual_\set{I}^z\right)$ can be upgraded to $\left(\dual_\set{I}\right)$ in which the sum of the reduced costs for \set{I} and the objective function of $\left(\dual\right)$ are gathered in a new objective function. The preliminary computation of $z^*$ is not necessary anymore.
		$$\left(\dual_{\set{I}}\right)\left\{
	\begin{array}{rrll@{\quad}l@{\qquad}l@{\qquad}l}
		\multicolumn{4}{l}{\max \omega_{\set{I}} = ub + \frac{1}{|\set{I}|}\sum\limits_{kl\in \set{I}}r_{kl}(u)}\\
		\text{s.t. } 	
		& uA & \le & c & & \\
		& u & \ge & 0 &\\
	\end{array}\right.$$
	
	
	\begin{property}[Usefulness of $\left(\dual_{\set{I}}\right)$]
		\label{prop: Usefulness of D_I} 
		
		Under Assumption \ref{feasAssump}, if $u^*$ is an optimal solution for $\left(\dual_{\set{I}}\right)$, we have		
		$\forall ij \in \set{I}$, 
		$$u^{*}b + r_{ij}(u^*) = z_{|ij}^*$$
		
	\end{property}
	\begin{proof}
		\prooflabel{Usefulness of D_I}
		We prove first that  $\omega^*_{\set{I}} = z^* +\frac{1}{|\set{I}|}\sum\limits_{kl\in\set{I}}R_{kl}$:
			\begin{itemize}
				\item 
				An optimal solution $u_1$ for $\dual_{\set{I}}^z$ is also feasible for $\dual_{\set{I}}$. By construction,  $u_1 b= z^*$ and $\forall kl \in \set{I}$, $r_{kl}(u_1) = R_{kl}$. 
				
				Therefore, under Assumption \ref{feasAssump},
				$\begin{aligned}[t]
					\omega_{\set{I}}^* &\ge u_1 b + \frac{1}{|\set{I}|} \sum\limits_{kl\in \set{I}} r_{kl}(u_1) \ge z^* + \frac{1}{|\set{I}|}\sum_{kl\in\set{I}} R_{kl}
				\end{aligned}$
				\item 
				A feasible solution $u_2$ for $\dual_{\set{I}}$ is also feasible for \dual.
				
				Thus, $\begin{aligned}[t]
					&& u_2 b +r_{kl}(u_2) &\le z^* + R_{kl}	\quad \forall kl \in \set{I}\\
					\implies && |\set{I}| \left(u_2 b\right) +\sum\limits_{kl\in\set{I}}r_{kl}(u_2) &\le |\set{I}|\,z^* + \sum\limits_{kl\in\set{I}}R_{kl}\\
					\implies && u_2 b + \frac{1}{|\set{I}|}\sum\limits_{kl\in\set{I}}r_{kl}(u_2) &\le z^* +\frac{1}{|\set{I}|} \sum\limits_{kl\in\set{I}}R_{kl}\\
					\implies && \omega_{\set{I}} & \le z^* +\frac{1}{|\set{I}|} \sum\limits_{kl\in\set{I}}R_{kl}\\
					\implies && \omega^*_{\set{I}} & \le z^* +\frac{1}{|\set{I}|} \sum\limits_{kl\in\set{I}}R_{kl}\\					
				\end{aligned}$
			\end{itemize}		
			Thus, $\omega^*_{\set{I}} = z^* +\frac{1}{|\set{I}|}\sum\limits_{kl\in\set{I}}R_{kl}$\\
			We can now prove the property. Let $u_3$ be an optimal solution for $\left(\dual_{\set{I}}\right)$,
			
			$\forall ij \in B$, $u_3 b + r_{ij}(u_3) \le z^*+R_{ij}$.
			
			\noindent Suppose that $\exists \widetilde{kl}\in \set{I}$ \st $u_3 b + r_{\widetilde{kl}}(u_3)< z^*+R_{\widetilde{kl}}$.
			
			\noindent$|\set{I}|\,\omega^*_{\set{I}} ={|\set{I}|} \, z^* +\sum\limits_{kl\in\set{I}}R_{kl}$, implies $\sum\limits_{kl\in \set{I}\setminus \widetilde{kl}}\left(u_3 b + r_{kl}(u_3) \right)>\sum\limits_{kl\in\set{I}\setminus \widetilde{kl}} \left( z^*+R_{kl}\right)$ which is impossible.
			
			Thus, $\forall kl \in \set{I}, u_3 b + r_{kl}(u_3) = z^* + R_{kl} = z_{|kl}^*$
	\end{proof}

	Moreover, the original optimal value $z^*$ is available as a side product when a set of incompatible pairs is known to contain at least one pair of an optimal solution.
	\begin{corollary}
		If \set{I} contains at least one pair of $B$ belonging to an optimal solution of $\left(\primal\right)$ and $u^*$ is an optimal solution for $\left(\dual_{\set{I}}\right)$, then
		
		$$z^* = u^*b + \min\limits_{ij\in\set{I}} r_{ij}(u^*)$$
	\end{corollary}

	Property \ref{prop: Usefulness of D_I} and its corollary gives a simple algorithm to compute a lower bound for $Z$ and to achieve arc-consistency:
	

	\begin{algorithm}[H]
		\caption{ACbyLP \label{algoACbyLP}}
		\begin{algorithmic}[1]
			
			\STATE Unmark all variable-value pairs $ij \in B$
			\STATE $Zlb = +\infty$
			\STATE let $\set{T}$ be a complete set of incompatible pairs ($\cup_{\set{I} \in \set{T}} \set{I} = B$)
				\FOR {each $\set{I} \in \set{T}$}
				\IF {$\set{I}$ has unmarked pairs}
					\STATE Compute $\tilde{u}$ an optimal solution of $\left(\dual_{\set{I}}\right)$ and $w= \tilde{u}b$
					\STATE $Zlb = \min \{ Zlb ~;~  w + \min\limits_{ij\in \set{I}} r_{ij}(\tilde{u})\}$
					\FOR {$kl\in B$, $kl$ unmarked}
						\STATE \textbf{if} $r_{kl}(\tilde{u}) + w > \overline{Z}$ \textbf{ then} mark $kl$ as inconsistent.
					\ENDFOR	
					\FOR {$kl\in \set{I}$, $kl$ unmarked}
						\STATE \textbf{if} $r_{kl}(\tilde{u}) + w \leq \overline{Z}$ \textbf{then } mark $kl$ as consistent.
					\ENDFOR	
				\ENDIF			
			\ENDFOR
			\STATE Update $\underline{Z}$ with $Zlb$.
		\end{algorithmic}				
	\end{algorithm}
	
	Algorithm \ref{algoACbyLP} considers the sets $\set{I}$ of  incompatible pairs, one by one. For each $\set{I}$, $\left(\dual_{\set{I}}\right)$ is solved to get the exact reduced cost of the pairs of $\set{I}$. Note that the dual solution obtained is used to filter the entire domains. The set of pairs whose status consistent/inconsistent have been definitely established are marked. The consistency of some additional pairs than $\set{I}$ might be establised so line 9 could be extended with additional conditions based, for instance, on Theorem \ref{ac slackness} or sensitivity analysis. Algorithm \ref{algoACbyLP} does not specify how the family of incompatible sets should be built but the set of domains can be used by default. Moreover, the algorithm can be stopped at any time providing valid filtering for the whole domains. The order to consider the sets of $\set{T}$ is also left unspecified and many strategies can be imagined. 
	
	We believe that an anytime algorithm is key for very costly global constraints where arc-consistency is rarely worth a high runtime complexity such as $O(n^3)$. See for instance the discussion in \citep{DBLP:journals/constraints/CauwelaertS17} where the  arc-consistency algorithm for \constraint{MinWAllDiff} is found too costly and the filtering of \citep{focacci1999cost} used as a baseline is too weak.  Reduced costs based filtering techniques could be a very good framework to design anytime and adaptive consistency algorithms \citep{balafrej2016adapting}.

	\section{Application to the maximum Weighted Independent Set (WIS) in chordal graph}
\label{section: MWIS}

The general results presented so far give a methodology to design a complete propagator for a weighted global constraint. Using LP, it can provide a way to efficiently prototype the arc-consistency algorithm and already leads to a non-trivial approach. However, it can also give insight on how to derive AC from a dedicated combinatorial algorithm that already relies on duality and complementary slackness as an optimality termination condition. Typical examples are Primal-Dual algorithms \citep{papadimitriou1998combinatorial}. We propose to demonstrate this in the present section. We apply our results to the maximum weighted independent set problem (WIS) in chordal graphs, since an ideal formulation and a combinatorial algorithm are known. The resulted algorithm is novel and improves over a straightforward AC algorithm \citep{apeloig:dumas-04936273}.

\paragraph{Preliminaries about chordal graphs} Let $G = (V, E)$ be an undirected graph defined over a set $V$ of vertices and a set $E$ of edges. Given a vertex $v \in V$, we denote by $N(v)\subseteq V$ its \textit{neighborhood} (that is, the set of vertices adjacent to $v$ in $G$).
We call a set of pairwise adjacent vertices a \textit{clique}, and a set of pairwise non adjacent vertices a \textit{stable set} (or \textit{independent set}).
A vertex is called \emph{simplicial} if its neighborhood is a clique. A \textit{perfect elimination ordering} (\textit{peo}) of $G$ is an ordering $\{v_1, \ldots, v_{|V|}\}$ of the vertices of $G$ such that for every $1 \le i \le |V|$, $v_i$ is simplicial in the graph induced by $\{v_i, v_{i + 1}, \ldots, v_{|V|}\}$ (or equivalently, $N(v_i) \cap \{v_{i + 1}, \ldots, v_{|V|}\}$ is a clique).
As an example the ordering $h,f,i,a,c,b,d,g,e$ is a \emph{peo} for the graph of Figure \ref{MWIS_AC} whereas $h,b,\ldots$ is not.

A graph is said to be \textit{chordal} if it contains no induced cycle of length at least four (in other words, any cycle of length at least four has a chord). Moreover, $G$ is chordal if and only if $G$ has a \textit{peo} \citep{Fulkerson1965IncidenceMA}.
Chordal graphs are known to be perfect graphs which have a chromatic number equal to the size of a maximum clique. Moreover,  a chordal graph admits at most $|V|$ inclusion-wise maximal cliques.
We refer the reader to \citep{ramirez2001perfect, introChordalGraph} for more details about chordal graphs.

\paragraph{A global constraint for weighted independent set} 
Let $G = (U, E)$ be a graph on $n$ vertices $\{v_1, \ldots, v_n\}$. Let $\mathcal{Q}$ be a set of cliques of $G$ covering all edges. More precisely, for every $e \in E$ there exists a clique $Q \in \mathcal{Q}$ such that $e$ is an edge of $Q$.
Let $\{X_1,\ldots, X_n\}$ be a set of 0/1 variables. For every $1 \le i \le n$, $X_i = 1$ means that the vertex $v_i$ is included in a stable set of $G$.
In addition, there is for every $1 \le i \le n$, a weight (or cost) $w_i \in \mathbb{N}$ associated to $v_i$.
In order to follow the previous framework, we could define the assignment cost $c$ as $c_{i0} = 0$ (the cost of assigning the value $0$ to the variable $X_i$) and $c_{i1} = w_i$ (the cost of assigning the value $1$ to the variable $X_i$) for every $1 \le i \le n$. We consider the constraint $\constraint{MaxWIS}(X_1,\ldots, X_n, Z, G, c)$ which enforces the variables $X_1,\ldots, X_n$ to define the characteristic vector of an independent set (IS) of $G$ of total weight above the cost variable $Z$.
However, since $X_i$ is a 0/1 variable, we have $c_{iX_i} = w_iX_i$ for every $1 \le i \le n$. Therefore, we will directly use $w$ instead of $c$ in the following. Note also that the problem is stated as a maximization problem since it is often considered and presented as such in the literature. $\constraint{MaxWIS}$ is equivalent to the following constraint network:
	$$\{X \in \{0,1\}^n : \sum\limits_{v_i \in Q} X_i \leq 1\:\: \forall Q \in \mathcal{Q}, \:\:\:\: \sum^{n}\limits_{i=1} w_{i}X_i \geq Z \}$$

\begin{example}[Arc-consistency for $\constraint{MaxWIS}$]

\slshape\sffamily~
\begin{figure}[H]
\centering
    \includegraphics[width=9cm]{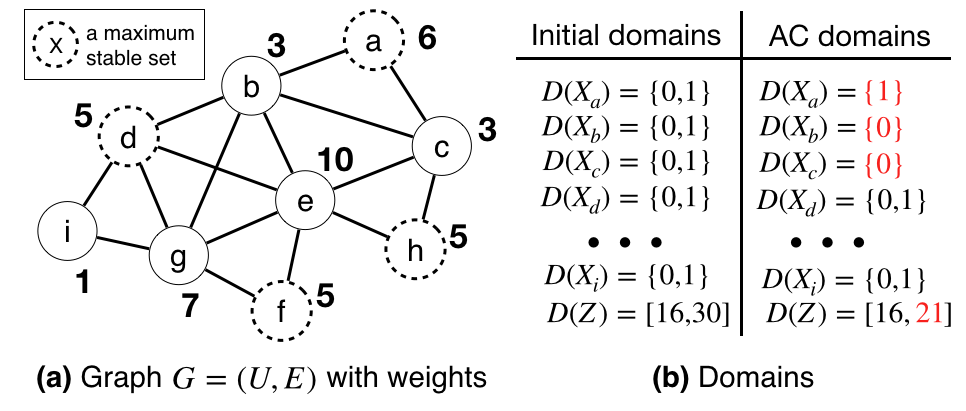}
      \caption{Example of arc-consistent domains for $\constraint{MaxWIS}(X_a,\ldots,X_i,Z,G,c)$ with an initial domain for the cost variable $D(Z) = [16,30]$.}
      \label{MWIS_AC}
\end{figure}
Figure \ref{MWIS_AC} gives an example of the filtering expected from arc-consistency for $\constraint{MaxWIS}$. A maximum weighted independent set has weight 21 in this example thus $\overline{Z} = 21$. Arc-consistent domains are shown on the right where vertices $b$ and $c$ have been instantiated to 0 (they do not belong to any stable stable set of weight above $\underline{Z} = 16$) and vertex $a$ is mandatory to reach 21 so $D(X_b) = D(X_c) = \{0\}$ and $D(X_a) = \{1\}$.
\end{example}
 
Two preliminary observations are required to state the general lines of the AC algorithm.\\
\indent Firstly, note that the arc-consistency closure  (ac-closure) of the clique inequalities alone, that is, $\sum_{v_i \in Q} X_i \leq 1,\:\:\: \forall Q \in \mathcal{Q}$, ensures that the network restricted to the cliques is globally consistent. A single inequality enforces all variables to 0 as soon as a variable within its scope is grounded to $1$ \emph{i.e} instantiated to value one. It only filters values 1 and only when a variable is grounded to 1. So, no fix point is needed since a value removal cannot trigger another one. Once this ac-closure is done, any remaining value of an ungrounded domain belongs to a feasible stable set (including or excluding the corresponding vertex since none of its neighbors have been selected yet). As a result, Assumption \ref{feasAssump} made in Section \ref{section: AC for global constraint with costs} that each remaining value belongs to at least one support (of the constraint without costs) holds after this process.\\
\indent Secondly, to deal with the weights, we can focus on the global consistency of values 1 alone. Focusing on values 1 consists in eliminating vertices that do not belong to any stable set of weight greater or equal to $\underline{Z}$. Once this is done, a vertex $v_i$ might be mandatory (filtering value 0)  only if all its neighbors have been forbidden. If not, the previous reasoning on forbidden vertices would not be complete. But if all its neighbors have been ruled out, it necessarily belongs to a maximum stable set and the loss for \textbf{not} including the vertex $v_i$ is readily available with $R_{i0} = -w_i$ and $z^*_{|i0}=z^* - w_i$. Thus if $z^* - w_i < \underline{Z}$, value 0 is filtered from $X_i$ and the knowledge of the exact reduced cost for values 1 alone is enough. Therefore, using the previous notations, we are considering a complete family of dual solutions (see Definition \ref{completeSetDEF}) with respect to $H = \{i1 | v_i \in U\}$ (pair $i1$ refers here to a pair $ij$ with $j=1$).

The AC algorithm can be summarized as follow:
\begin{enumerate}
\item Perform the ac-closure of the clique inequalities.
\item Compute the value $z^*$ of a maximum weighted stable set as well as the exact reduced costs of values 1. Update $\overline{Z}$ to $z^*$. Finally, for all $1 \le i \leq n$, if $z^* + R_{i1} < \underline{Z}$, remove value 1 from $D(X_i)$. 
\item For any ungrounded vertex $v_i$ ($D(X_i) = \{0,1\}$) whose neighbors have been grounded to 0 ($D(X_j) = \{0\},\:\:\forall v_j \in N(v_i)$) and such that $z^* - w_i < \underline{Z}$ ($R_{i0} = - w_i$), remove value 0 from $D(X_i)$.
\end{enumerate} Note that the last step alone is setting variables to 1 but only when all neighbors have been already decided so that there  is no need of a fix-point and the presented approach is idempotent.
Step 1 takes $O(|\mathcal{Q}|n)$ with a very simple approach and step 3 is in $O(n)$. Step 2 is the key step where we intend to demonstrate the usefulness of the previous general results. To do so, we solve this problem for chordal graphs. We now focus on Step 2 and propose two ways to implement it: using a general LP (Section \ref{sectionLPstable}) and using a dedicated algorithm  (Section \ref{sectionFrankstable}).

For the remaining of this section we consider that $G$ is a chordal graph, that $\{v_1, \ldots, v_n\}$ is a \emph{peo} of $G$ and that $\mathcal{Q}$ is the set of (inclusion-wise) maximal cliques of $G$. It is known that $|\mathcal{Q}| \leq n$ for chordal graphs.

\subsection{Application of the methodology using LP}
\label{sectionLPstable}

Consider the following primal and dual formulations (respectively 
$(\primal_{\constraint{WIS}})$ and $(\dual_{\constraint{WIS}})$) of the maximum weighted independent set problem:\\

\noindent$\begin{array}{l|l}	
(\primal_{\constraint{WIS}}) \qquad\qquad\qquad\qquad\qquad\qquad\:\:\:\:\:  & 	 (\dual_{\constraint{WIS}}) \\
\left\{				
	\begin{array}{rrll@{\quad}l@{\quad}l@{\quad}l}
	\multicolumn{6}{l}{\max z = \sum\limits_{v_i\in U}  w_ix_i}\\
		\text{s.t. } 	
		&\sum\limits_{v_i \in Q} x_i& \le & 1 & \forall Q \in \mathcal{Q}\\
		& x_{i} & \geq &  0 & \forall v_i \in U\\
\end{array}\right.
&
\left\{				
			\begin{array}{rrll@{\quad}l@{\qquad}l@{\qquad}l}
				\multicolumn{6}{l}{\min w = \sum\limits_{Q \in \mathcal{Q}} u_Q}\\
				\text{s.t. } 	
				& \sum\limits_{Q \in \mathcal{Q}|v_i \in Q} u_Q& \ge & w_i & \forall v_i \in U\\
				& u_Q & \ge & 0& \forall Q \in \mathcal{Q}
\end{array}\right.
\end{array}
$\\
\newline

\indent The polytope $P = \{x \geq 0 : \sum\limits_{v_i \in Q} x_i \le 1,  \forall Q \in \mathcal{Q}\}$ is known to be the stable set polytope (the convex hull of stable sets of $G$) for perfect graphs (Fact 9.4 in \cite{ramirez2001perfect}). The formulation $(\primal_{\constraint{WIS}})$ based on the clique inequalities of a perfect graph is an ideal formulation. A support of maximal cost for $\constraint{MaxWIS}$ can therefore be found by solving $(\primal_{\constraint{WIS}})$ with the simplex algorithm and a naive filtering algorithm can solve the formulation $n$ times for each $X_i=1$ (each $i1$). An improved algorithm is readily available using the previous results of this paper. It is enough to notice that each clique $Q \in \mathcal{Q}$ gives a set $\mathcal{I}_i = \{ i1 | \forall v_i \in Q\}$ of incompatible pairs where all exact reduced costs can be obtained in the same dual solution according to Corollary \ref{prop: incompatible edges}. Algorithm \ref{algoACbyLP} can therefore be run with a complete family $\mathcal{T}=\{\mathcal{I}_1,\ldots, \mathcal{I}_{|\mathcal{Q}|}\}$. Since $|\mathcal{Q}| \leq n$ for chordal graphs (and can be significantly smaller than $n$ in practice), this is an improvement. Recall that the value $z^*$ needed to update $\overline{Z}$ is available at the end of Algorithm \ref{algoACbyLP} (since it is a maximization problem in the present case).

\begin{property}[AC for max weighted independent set with $|\mathcal{Q}|$ LPs]
\label{stableLPAC}
Arc-consistency of \constraint{MaxWIS} can be achieved with the resolution of $|\mathcal{Q}|$ linear programs (\emph{i.e} solving $(\dual_{\set{I}})$ for each clique of $G$ at step 2 of the outline of the filtering algorithm).
\end{property}


\subsection{Application of the methodology using a dedicated combinatorial algorithm}
\label{sectionFrankstable}

A maximum weighted independent set in a chordal graph can be computed in linear time since the seminal work of \citep{PolyMWISFrank}. This algorithm relies on duality theory and our results suggest to investigate whether we can adapt Frank's approach to design a dedicated algorithm for ($D_{\set{I}}$). We will show how the exact reduced costs of the vertices of a clique can be found with one call to Frank's algorithm so that AC can be achieved in $O(|\mathcal{Q}||E|)$ time which improves over a naive $O(n|E|)$ since $|\mathcal{Q}| \leq n$. Let's start with Frank's algorithm.

\paragraph{Frank's algorithm} The algorithm considers a \emph{peo} of the chordal graph and runs in two phases: forward and backward. Let us explain how the algorithm works using the \emph{peo} $\{v_1, \ldots, v_n\}$. The weight of each vertex is updated during the course of the algorithm and we refer to $w'_i$ as its current weight to distinguish it from its original weight $w_i$. Initially, $w'_i$ is set to $w_i$.
\begin{itemize}


\item \textit{Forward phase:} The algorithm iterates over the vertices in the order of the \emph{peo}. Whenever it reaches a vertex $v_k$ with $w'_k > 0$, it marks it \textbf{in red}. Then, it sets $w'_l$ to $w'_l - w'_k$ for every $l > k$ such that $v_kv_l \in E$.
Note that the original Franck's algorithm sets $w'_k$ to $0$ when $v_k$ is marked in red and $w'_l$ to $0$ as soon as it becomes negative, but we shall keep the values which will be useful later on. The set of vertices marked in red at the end of the forward pass is denoted by $\mathcal{R}$. Moreover, we will denote by $\mathcal{Q}_\mathcal{R}$ the set of the corresponding cliques in the \textit{peo}, that is, $\mathcal{Q}_\mathcal{R} = \{N(v_k)\cap \{v_{k+1},\dots,v_n\} | v_k \in \mathcal{R}\}$. Keep in mind that each clique is related to a vertex $v_k$ at position $k$ in the \emph{peo} and weight $w'_k$ at the end of Frank's algorithm.

\item \textit{Backward phase:}  The algorithm, goes through all the vertices marked in red in the reverse order and marks in \textbf{blue} each vertex that is not adjacent to a vertex already marked in blue. Let $S$ be the independent set of vertices marked in blue obtained after this phase. 
\end{itemize}
\begin{example}[Frank's algorithm]
\slshape\sffamily~
\begin{figure}[H]
    \centering
    \includegraphics[width=11cm]{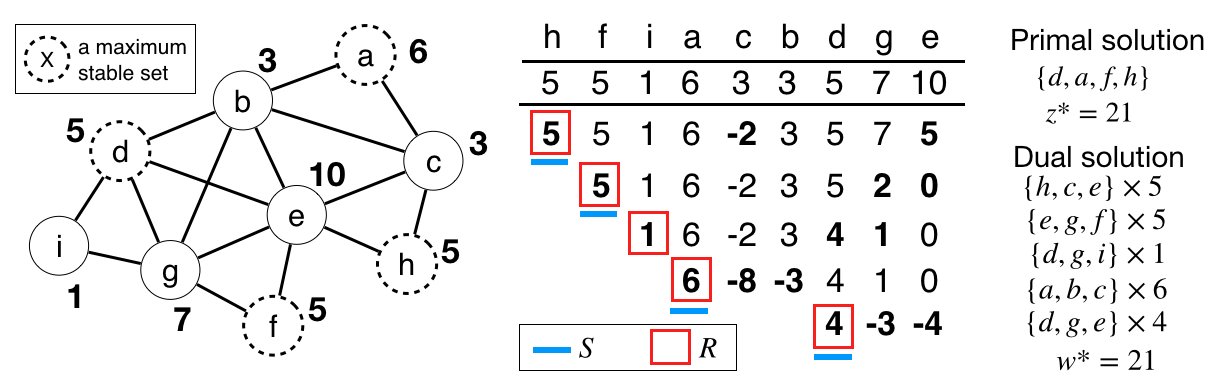}
      \caption{Example with a weighted graph (left) and the execution of Frank's algorithm using the \emph{peo} $\{h,f,i,a,c,b,d,g,e\}$. An optimal stable set $S = \{d,a,f,h\}$ is found.}
    \label{fig:frankrun}
\end{figure}

Figure \ref{fig:frankrun} shows an execution of Frank's algorithm. Since $w'_h=5>0$, we mark $h$ in red and for each adjacent vertex of $h$ among $\{f,\dots,e\}$ we remove $w'_h$. The corresponding \textit{clique} is $\{h,c,e\}$ so that $w'_c$ and $w'_e$ become respectively -2 and 5. At the next iteration, $f$ is selected since $w'_f=5>0$, marked red and $w'_g$ and $w'_e$ become respectively 2 and 0. And so on. At iteration 5, vertices $c$ and $b$ are skipped ($w'_c=-8\leq0$, $w'_b=-3 \leq0$) and $d$ is marked red, $w'_g$ becomes -3 and $w'_e$ becomes -4. Once all vertices have been considered, we have $\mathcal{R}=\{h,f,i,a,d\}$.  We go through the red vertices in the reverse order: $d$ is selected, $a$ is selected (not adjacent to $d$), $i$ is discarded because it is adjacent to $d$, $f$ and $h$ are selected. We obtain the independent set $S=\{d,a,f,h\}$ for a total weight of $z^*=21$.

We can observe that $\mathcal{Q}_\mathcal{R}=\{\{h,c,e\}, \{e,g,f\},\{d,g,i\},\{a,b,c\},\{d,g,e\}\}$ gives a feasible dual solution of $(\dual_{\constraint{WIS}}):u_{\{h,c,e\}} = 5$, $u_{\{e,g,f\}} = 5$, $u_{\{d,g,i\}} = 1$, ${u_{\{a,b,c\}} = 6}$, $u_{\{d,g,e\}} = 4$. A clique $Q \in \mathcal{Q}_\mathcal{R}$ is given the weight $w'_i$ of the vertex $v_i$ of the corresponding iteration, in other words $u_Q = w'_i$. Feasibility of the dual solution $u$ can be checked. For instance for vertex $e$, we have $u_{\{h,c,e\}} + u_{\{e,g,f\}} + u_{\{d,g,e\}} = 5+5+4 \geq w_e = 10$. Note also that $\sum_{Q \in\mathcal{Q}_R} u_Q = 21$ which prove the optimality of $S$ by strong duality.
\end{example}

Since $G$ is chordal, the time complexity of finding a \emph{peo} of $G$ is $O(|E| + |U|)$ 
\citep{tarjanEtAL1984} and the forward/backward phase runs in linear time $O(|E|)$.
The optimality of Frank's algorithm can be established by showing that it ends with a primal feasible  solution of $(\primal_{\constraint{WIS}})$ and a dual feasible solution of $(\dual_{\constraint{WIS}})$ which have the same objective value. 
Interestingly, the algorithm explicitly outputs the reduced costs of each vertex corresponding to the optimal dual solution $u^*$ found. The reduced cost of each vertex $v_i$ (variable $x_i$ of $(\primal_{\constraint{WIS}})$) with $u^*$ is given by:
$$r_{i1}(u^*) = w_i-\sum_{Q\in \mathcal{Q}|v_i\in Q} u^*_{Q}$$

\begin{example}[Frank's algorithm (continued)]
 In the example of Figure \ref{fig:frankrun}, we have $r_{i1}(u^*)=0$, $r_{c1}(u^*)=-8$, $r_{b1}(u^*)=-3$, $r_{g1}(u^*)=-3$, $r_{e1}(u^*)=-4$. Reduced costs for vertices $c,g,e$ are actually exact. Typically, it is not exact for $b$ since the best possible independent set including $b$ has a total weight of 14 so that $R_{b1}=-7$.
 \end{example}
 
We can show that the reduced cost of the vertices located after the last red vertex in the ordering are guaranteed to be exact at the end of the algorithm (Property \ref{endClique}). The key idea is then to build a \emph{peo} ensuring that a maximal clique ends the ordering, starting with its vertex of highest weight (Property \ref{obsCliqueEnd}). If that is possible, Frank's algorithm will provide the exact reduced cost of the corresponding maximal clique and arc-consistency can be established with $O(|\mathcal{Q}|)$ calls to Frank's algorithm (Theorem \ref{theoremACqm}). In other words, we have a dedicated algorithm to solve the generic problem $(D_{\set{I}})$ of Section~\ref{section: LP base algo}. The remaining of the section is dedicated to proving each step of the algorithm. We denote by $N_\mathcal{R}(v_k) = N(v_k) \cap \mathcal{R}$ the set of vertices marked in red and adjacent to $v_k$. Similarly, we denote by $N^-_\mathcal{R}(v_k) = N_\mathcal{R}(v_k) \cap \{v_1, \dots, v_{k-1}\}$ the set of vertices marked in red, adjacent to $v_k$ and preceding it in the \textit{peo}.

\begin{lemma}[Weight of vertices after the last red]
\label{lemMargCost}
After an execution of Frank's algorithm with a \emph{peo} $\{v_1, \ldots, v_n\}$, we have for every vertex $v_k$ after the last red vertex: 
$$z_{|k1}^*=\sum\limits_{v_i\in \mathcal{R}\backslash N_\mathcal{R}(v_k)}w_i'+w_k$$
\end{lemma}
\begin{proof}
Let $v_k$ be a vertex which is after the last red vertex in the \emph{peo} after an execution of Frank's algorithm. One can observe that since every vertex in $\mathcal{R}$ precedes $v_k$ in the \emph{peo}, we have $N_\mathcal{R}^-(v_k)=N_\mathcal{R}(v_k)$. Now let us execute Frank's algorithm on the graph induced by $U \setminus N(v_k)$ using the same $\emph{peo}$ (minus the vertices in $N(v_k)$). We denote by $\mathcal{R}^*$ the red vertices and by $w''_i$ the final weight of $v_i$ for every $v_i$ in $U \setminus N(v_k)$ in this execution.
Let us show that $w''_k = w_k$ and that $w''_i = w'_i$ for every $v_i \in U' := U \setminus (N(v_k) \cup \{v_k\})$. Since $v_k$ does not have any neighbor in $U'$, we have $w''_k = w_k$. 

If $N^-_\mathcal{R}(v_i) = N^-_{\mathcal{R}^*}(v_i)$ for every $v_i \in U'$, then since all of the vertices of $U'$ are non adjacent to $v_k$, we have $w''_i = w'_i$ for every $v_i \in U'$. Suppose that it is not the case and let $v_j$ be the first vertex (in the \emph{peo}) of $U'$ such that $N^-_\mathcal{R}(v_j) \neq N^-_{\mathcal{R}^*}(v_j)$. In order for $v_j$ to be the first vertex of this sort, it must be adjacent to some neighbor $v_l$ of $v_k$, such that $l < j$. But then, by definition of a \emph{peo}, since $v_j$ and $v_k$ are both right-neighbors of $v_l$, $v_j$ must be adjacent to $v_k$ which is a contradiction. Since the value $z^*_{|k1}$ is exactly the cost of a solution given by Frank's algorithm on the graph induced by $G$ on $U \setminus N(v_k)$ the proof of the lemma is complete.
\end{proof}

\begin{property}[Exact reduced costs]\label{endClique}
Consider a chordal graph $G$ and an execution of Frank's algorithm over a given \textit{peo}. Reduced costs of all vertices after the last red vertex in the \emph{peo} are exact.
\end{property}
\begin{proof}
    Assume Frank's algorithm has been executed and consider a vertex $v_k$ after the last red with a reduced cost of $w'_k$. Since $v_k$ is after the last red vertex in the \emph{peo} we can apply Lemma~\ref{lemMargCost} to get $z^*_{|k1} = \sum\limits_{v_i\in \mathcal{R}\setminus N_\mathcal{R}(v_k)} w'_i + w_k$. Furthermore, we have  $w'_k = w_k - \sum\limits_{v_i \in N_\mathcal{R}(v_k)} w'_i$. This gives us
    \begin{align*}
        z^*_{|k1}   &= \sum\limits_{v_i\in \mathcal{R}\setminus N_\mathcal{R}(v_k)} w'_i + \sum\limits_{v_i \in N_\mathcal{R}(v_k)} w'_i + w'_k \\
                    &= \sum\limits_{v_i \in \mathcal{R}} w'_i + w'_k \\
                    &= z^* + w'_k.
    \end{align*}
    Hence, $w'_k$ is the exact reduced cost of $x_k$.
\end{proof}


It turns out that exact reduced costs of additional vertices (other than the ones after the last red) can be computed and the criterion of Property \ref{endClique} can be generalized. The reduced cost of a vertex $v_k$ whose red neighborhood occurs before position $k$ in the \emph{peo} (\emph{i.e} $N_\mathcal{R}(v_k) \cap \{v_{k+1},\ldots, v_n\} = \emptyset$) is exact. Although this would be a stronger point, it is not needed to complete our result. We turn our attention to show that we can build a family of \emph{peos} such that each ends with a clique of $G$.

\begin{property}\label{obsCliqueEnd}
Given a chordal graph $G$, for any clique $Q$ of $G$, there exists a \emph{peo} that ends with the vertices of $Q$.
\end{property}
\begin{proof}
It is known that every chordal graph G which is not a clique contains at least two non-adjacent \textit{simplicial} vertices \citep{Dirac61}. A given clique can therefore always be postponed at the very end. During the process, when a \textit{simplicial} vertex must be selected, either the resulting induced graph is $Q$ or it exists $v \notin Q$ which is a \textit{simplicial} vertex and can be selected first.
\end{proof}

Note that the order of the vertices in an ending clique $Q$ does not matter (any order of these vertices gives a valid \emph{peo}). We will therefore sort the vertices of the ending clique during the execution of Frank's algorithm, and more precisely right before marking the first vertex of the clique. The vertex of the clique $Q$ that will be put in first position is $v_k=$ arg$\max_{v_i\in Q} w'_i$, where $w'$ denotes the current weights of the vertices during the execution. This operation is made only once during the execution of Frank's algorithm and is made in time $O(|Q|)$. Thus, the complexity of the algorithm is unchanged.


\begin{theorem}[Complexity of arc-consistency with reduced costs]\label{theoremACqm}
    Arc-consistency can be achieved in $O(|\mathcal{Q}||E|)$.
\end{theorem}

\begin{proof}
According to Property~\ref{obsCliqueEnd}, we can create for each clique $Q \in \mathcal{Q}$, a \textit{peo} such that all vertices of $Q$ are at the end.
During the execution, once all vertices but the ones of $Q$ have been considered, either all vertices of $Q$ have a non positive value or it exists at least one vertex $v_k\in Q$ such that $w'_k>0$. If such vertex exists, mark in red the vertex $v_k$ defined as arg$\max_{v_i\in Q} w'_i$ and compute the resulting weights of the others. In both cases, according to Property~\ref{endClique}, for each vertex of $Q$, we obtain its \textit{exact reduced cost}. Note that in the case where $w'_k>0$, $v_k$ will be the last red marked vertex. Thus, it will be marked in blue and will appear in at least one optimal solution, its exact reduced cost will be 0 ($R_{k1}=0$). Each vertex of $U$ appears in at least one clique of $\mathcal{Q}$. All exact reduced costs are obtained by running the algorithm once for each $Q \in \mathcal{Q}$ with such a \textit{peo}. Since building a peo and running Frank's algorithm are done in a linear time complexity, all exact reduced cost are computed in $O(|\mathcal{Q}||E|)$ (Step 2) and arc-consistency takes $O(|\mathcal{Q}||E|)$.
\end{proof}

The AC algorithm is illustrated on our running example in Annex (Example \ref{example_ac_algo}) where three \emph{peo} are used to achieve arc-consistency.

	\section*{Conclusion}
	We show that arc-consistency can be done for global constraints with assignment costs by solving $n$ linear programs in the worst case, one for each variable and that this bound is sharp. To our knowledge, it provides the first analysis of reduced cost filtering which has often been used in the past in CP starting with the work of \citep{focacci1999cost}.
	
	
	This analysis established a number of basic results relating reduced costs and AC by answering the following questions: does there always exist a dual solution that can assert the consistency of a value (property \ref{prop: exists u st r_kl=R_kl}) ? Given a dual solution, how do we know which values are proved consistent/inconsistent (Theorem \ref{ac slackness}) ? Can we identify simple sufficient conditions for a family of dual solutions to ensure arc-consistency (Corollary \ref{prop: incompatible edges}) ?
	We also showed how these ideas can be directly applied to design an AC algorithm for stable sets in chordal graphs which have numerous applications. 

	Two key results of this paper are probably the characterization given by Theorem \ref{ac slackness} which states complementary slackness conditions for exact filtering (as opposed to just optimality) and Corollary \ref{prop: incompatible edges} which immediately lead to practical algorithms as demonstrated on the weighted independent set problems on chordal graphs. 
	
	The proposed analysis assumes an ideal ILP formulation to ensure AC of integer problems. But the results open the way to the analysis of any linear relaxation to obtain the best possible reduced costs related to a given relaxation.
	We believe this analysis contributes to the general question on how to search the dual space to perform filtering \citep{sellmann2004theoretical}.
	    
	\section*{Acknowledgment}
We thank the reviewers for their encouraging comments and help, without which we would not have been able to complete this work.

\clearpage
\section*{Annex}

\noindent \textbf{Property 1.} For any dual optimal solution $u^*$ of $(\dual)$ and any pair $kl \in B$, we have $$0\le r_{kl}(u^*)\le R_{kl}$$
	\begin{proof}[Proof of Property \ref{prop: bounds of RC}]
Solution $u^*$ is feasible for $(\dual)$ so the slack of the dual constraint related to variable $x_{kl}$ must be positive or null: $r_{kl}(u^*) \geq 0$. To show the second ($r_{kl}(u^*) \leq R_{kl}$), let $\tilde{x}^*$ be an optimal solution of $\left(\primal_{|kl}\right)$, the problem restricted with $X_k = l$. Solution $\tilde{x}^*$ exists by Assumption \ref{feasAssump}. By strong duality of $(\primal)$ and $(\dual)$:
	    
	    \indent$\begin{aligned}[t]
					z^* 
						&= u^*b\\
					z^*	& \le u^*(A\tilde{x}^*) \qquad && (\textrm{$\tilde{x}^*$ is a feasible solution of $(\primal$)})\\							
					z^*	& \le (u^*A)\tilde{x}^* && \\		
					z^*	&\le  (c - r(u^*))\tilde{x}^* && \left(\text{from the definition of the dual slack $r(u) = c - uA$}\right)\\
					z^*	&\le z^*_{|kl} - r(u^*)\tilde{x}^*&& (c\tilde{x}^*=z^*_{|kl} \text{ since $\tilde{x}^*$ is optimal for $\left(\primal_{|kl}\right)$})\\
	    \end{aligned}$\\
	
					\noindent Recall that $R_{kl} = z^*_{|kl} - z^*$. Thus from the last inequality, we have :\\
				\indent$\begin{aligned}[t]
						&R_{kl} \ge r(u^*)\tilde{x}^*\\
						&R_{kl}  \ge r_{kl}(u^*) & (\text{since $\tilde{x}^*_{kl}=1$ and $\forall ij \in B$, $r_{ij}(u^*)\ge 0$ and $\tilde{x}^*_{ij}\ge 0$})
						\end{aligned}$
		\end{proof}

\begin{example}[Arc-consistency algorithm for $\constraint{MaxWIS}$]
\slshape\sffamily~
\begin{figure}[H]
\centering
    \includegraphics[width=11cm]{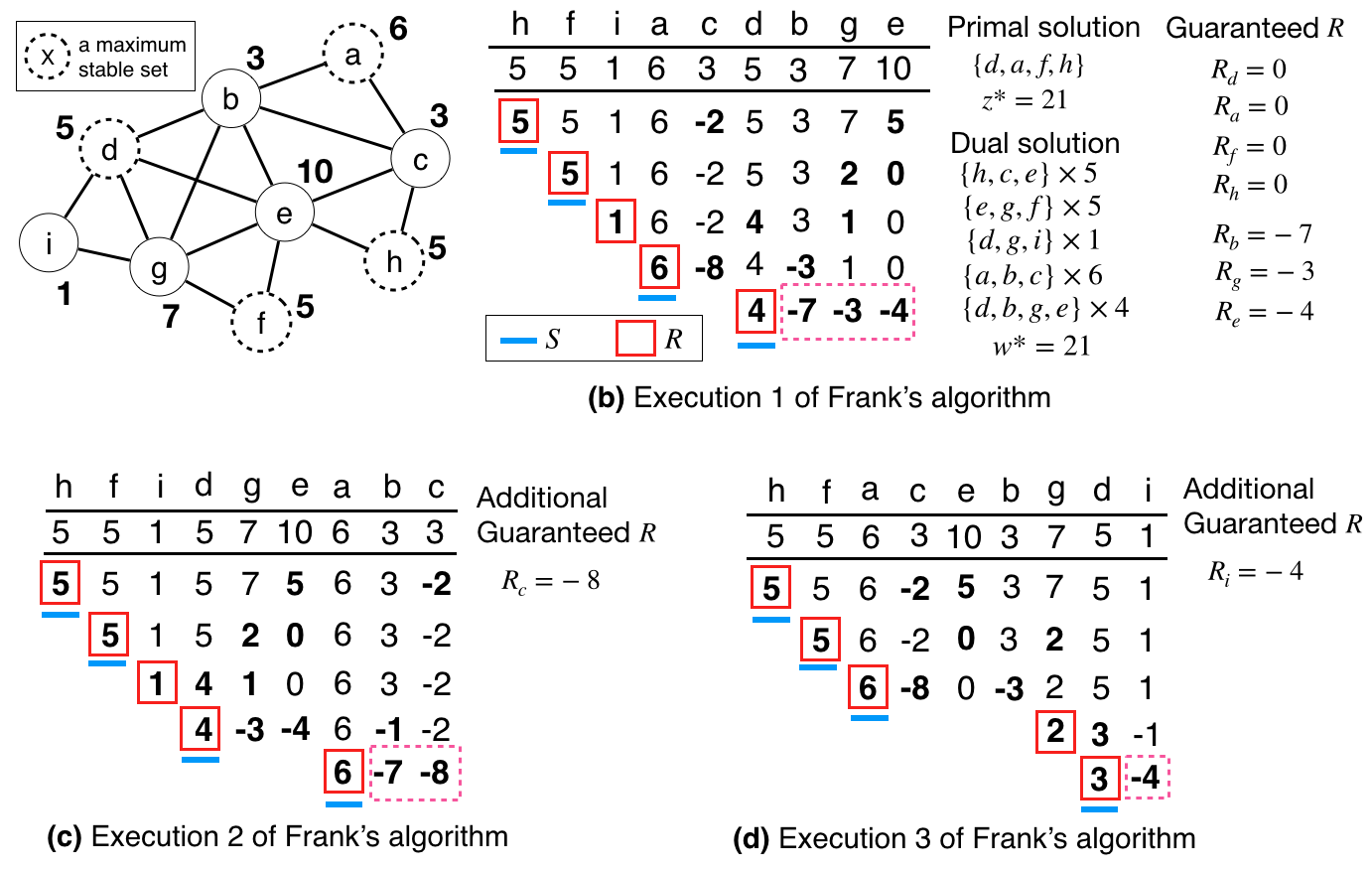}
      \caption{Example of the AC algorithm where all exact reduced costs are gathered after running Frank's algorithm for three \emph{peo}.}
      \label{MWIS_AC_Algoi}
\end{figure}
The first \emph{peo} targets the clique $\{d,b,g,e\}$, the second targets $\{a,b,c\}$ and the last puts $\{g,d,i\}$ at the end of the ordering. There is no need to continue since, at this stage all exact reduced cost have been obtained. Execution \textbf{(b)} is the same as in Figure \ref{fig:frankrun} where $d$ and $b$ have been swapped. During the execution, when all vertices of $U\setminus Q$ have been considered, arg$\max_{v_i\in Q} w'_i=v_d$. As $w'_d>0$, it is marked in red and $w'_b=w'_b-w'_d$, $w'_g=w'_g-w'_d$ and $w'_e=w'_e-w'_d$. Note that in execution \textbf{(d)}, $v_d$ should be taken before $v_g$. However, as only $R_i$ is needed, there is no difference in this case.

\label{example_ac_algo}
\end{example}

\bibliographystyle{apalike}
\bibliography{bibliography}

\end{document}